\documentclass{amsart}
\usepackage{amssymb}

\theoremstyle{definition}

\theoremstyle{remark}

\numberwithin{equation}{section}

\input{tcilatex}

\begin{document}
\title{Operator-valued Fourier multipliers in Besov spaces and its
applications}
\author{Veli Shakhmurov}
\address{Department of Mathematics, Okan University, Akfirat Beldesi, Tuzla,
34959, Istanbul, Turkey.}
\email{veli.sahmurov@okan.edu.tr}
\author{Rishad Shahmurov}
\address{Department of Mathematics, Yeditepe University, Kayishdagi Caddesi,
34755 Kayishdagi, Istanbul, Turkey}
\email{shahmurov@hotmail.com}
\subjclass[2000]{Primary 42B15, 46E40, 34G10 }
\date{}
\keywords{Fourier type spaces, Banach--valued Besov spaces, operator--valued
multipliers, interpolation of Banach spaces, boundary value problems,
differential--operator equations, operator convolution equations.}

\begin{abstract}
The present paper, is devoted to investigation of operator--valued Fourier
multiplier theorems from $B_{q_{1},r}^{s}$ to $B_{q_{2},r}^{s}$, optimal
embedding of Besov spaces, the separability and positivity of differential
operators. Here, we show that these differential operators generate analytic
semigroup.
\end{abstract}

\maketitle




\section*{1. Introduction, notations and background}

In recent years, Fourier multiplier theorems in vector--valued function
spaces have found many applications in embedding theorems of abstract
function spaces and in theory of differential operator equations, especially
in maximal regularity of parabolic and elliptic differential--operator
equations. Operator--valued multiplier theorems in Banach--valued function
spaces have been discussed extensively in [3,8,12,15, 17]. Boundary value
problems (BVPs) for differential--operator equations (DOEs) in $H$--valued
(Hilbert valued) function spaces and parabolic type convolution operator
equations (COEs) with bounded operator coefficient have been studied in
[1,2,6,7,9,13,14], and $\left[ \text{3}\right] $ respectively.

Let $E$ be a Banach space$,~x=(x_{1},x_{2},\cdots ,x_{n})\in \Omega \subset
R^{n}.$ $L_{r}(\Omega ;E)$ denotes the space of all strongly measurable $E$%
--valued functions that are defined on the measurable subset $\Omega \subset
R^{n}$ with the norm 
\begin{equation*}
\begin{array}{lll}
\Vert f\Vert _{L_{r}(\Omega ;E)} & = & \displaystyle\left( \int \Vert
f(x)\Vert _{E}^{r}dx\right) ^{\frac{1}{r}},~~~1\leq r<\infty , \\ 
\vspace{-3mm} &  &  \\ 
\Vert f\Vert _{L_{\infty }(\Omega ;E)} & = & \displaystyle\mbox{ess\,sup}%
_{x\in \Omega }[\Vert f(x)\Vert _{E}].%
\end{array}%
\end{equation*}

Let $S=S\left( R^{n};E\right) $ denote a Schwartz class i.e. the space of $E$%
-valued rapidly decreasing smooth functions on $R^{n},$ equipped with its
usual topology generated by semi-norms. Let $S^{\shortmid }\left(
R^{n};E\right) $ denote the space of all\ continuous linear operators $%
L:S\rightarrow E,$ equipped with the bounded convergence topology. Recall $%
S\left( R^{n};E\right) $ is norm dense in $B_{p,r}^{s}\left( R^{n};E\right) $
when $1\leq p,$ $r<\infty .$

Let $\alpha =(\alpha _{1},\alpha _{2},\cdots ,\alpha _{n}),$ where $\alpha
_{i}$ are integers. An $E$--valued generalized function $D^{\alpha }f$ is
called a generalized derivative in the sense of Schwartz distributions, if
the equality 
\begin{equation*}
<D^{\alpha }f,\varphi >~=~(-1)^{|\alpha |}<f,D^{\alpha }\varphi >
\end{equation*}%
holds for all $\varphi \in S.$

\vspace{3mm}The Fourier transform $F:S(X)\rightarrow S(X)$ is defined by 
\begin{equation*}
(Ff)(t)~\equiv ~\hat{f}(t)~=~\int\limits_{R^{N}}\exp (-its)f(s)ds
\end{equation*}%
is an isomorphism whose inverse is given by 
\begin{equation*}
(F^{-1}f)(t)~\equiv ~\check{f}(t)~=~(2\pi )^{-N}\int\limits_{R^{N}}\exp
(its)f(s)ds,
\end{equation*}%
where $f\in S(X)$ and $t\in R^{N}.$

It is known that 
\begin{equation*}
F(D_x^\alpha f)~=~(i\xi_1)^{\alpha_1}\cdots(i\xi_n)^{\alpha_n}\hat{f},~~~
D_\xi^\alpha(F(f))~=~F[(-ix_n)^{\alpha_1}\cdots(-ix_n)^{\alpha_n}f]
\end{equation*}
for all $f\in S^{\dagger}(R^n;E).$

\vspace{3mm}

Let $\mathbf{C}$ be the set of complex numbers and 
\begin{equation*}
S_{\varphi }~=~\{\lambda ;~\lambda \in \mathbf{C},~|\arg \lambda |\leq
\varphi \}\cup \{0\},~~~0\leq \varphi <\pi .
\end{equation*}%
A linear operator $A$ is said to be a $\varphi $--positive in a Banach space 
$E$, if $D(A)$ dense in $E$, and 
\begin{equation*}
\left\Vert (A+\lambda I)^{-1}\right\Vert _{B(E)}~\leq ~M(1+|\lambda |)^{-1}
\end{equation*}%
with $M>0,~\lambda \in S_{\varphi },~\varphi \in \lbrack 0,\pi )$; here $I$
is the identity operator in $E,~B(E)$ is the space of all bounded linear
operators in $E.$ Sometimes instead of $A+\lambda I$, we will write $%
A+\lambda $ and denote it by $A_{\lambda }.$ Let $E\left( A^{\theta }\right) 
$ denote the space $D\left( A^{\theta }\right) $ with graphical norm 
\begin{equation*}
\left\Vert u\right\Vert _{E\left( A^{\theta }\right) }=\left( \left\Vert
u\right\Vert ^{p}+\left\Vert A^{\theta }u\right\Vert ^{p}\right) ^{\frac{1}{p%
}},1\leq p<\infty ,-\infty <\theta <\infty .
\end{equation*}

\vspace{3mm}

\textbf{Definition 1}$.$ Let $A=A\left( t\right) $ be a uniformly positive
operator in $E,$ $u\in E(A)$ and $Au\in S\left( R;E\right) .$ Then, the
Fourier transformation of $A\left( t\right) $ is defined by 
\begin{equation*}
((FA)u)(\xi )=(\hat{A}u)(\xi )=(2\pi )^{-\frac{1}{2}}\int\limits_{R}e^{-it%
\xi }A(t)u\text{ }dt.
\end{equation*}

\vspace{3mm}

\textbf{Definition 2}$.$ Let $A=A\left( t\right) $ be a uniformly positive
operator in $E.$ Then, it is differentiable if for all $u\in $ $E\left(
A\right) $

\begin{equation*}
\left( \frac{d}{dt}A\right) u=A^{\prime }(t)u=\lim_{h\rightarrow 0}\frac{%
\left\Vert A(t+h)u-A(t)u\right\Vert _{E}}{h}<\infty .
\end{equation*}

\vspace{3mm}

\textbf{Definition 3}$.$ Let $A=A\left( t\right) $ be a uniformly positive
operator in $E$ and $u\in S\left( R;E(A)\right) $ and 
\begin{equation*}
(A\ast u)(t)=\int\limits_{R}A(t-y)u(y)\text{ }dy.
\end{equation*}

\vspace{3mm}

Let $y\in R,~m\in N$ and $e_{i},~i=1,2,\cdots ,n$ be standard unit vectors
of $R^{n}.$ Let 
\begin{equation*}
\begin{array}{l}
\Delta _{i}(y)f(x)~=~f(x+ye_{i})-f(x),~ \\ 
\vspace{-3mm} \\ 
\displaystyle\Delta _{i}^{m}(y)f(x)~=~\Delta _{i}(y)\left[ \Delta
_{i}^{m-1}(y)f(x)\right] ~=~\sum%
\limits_{k=0}^{m}(-1)^{m+k}C_{m}^{k}f(x+kye_{i}).%
\end{array}%
\end{equation*}%
Let 
\begin{equation*}
\Delta _{i}(\Omega ,y)~=~\left\{ 
\begin{array}{ll}
\Delta _{i}(y)f(x), & \mbox{for}~[x,x+mye_{i}]\subset \Omega \\ 
\vspace{-3mm} &  \\ 
0, & \mbox{for}~[x,x+mye_{i}]\notin \Omega .%
\end{array}%
\right.
\end{equation*}%
Let $m$ be integer, $s$ be positive number, and 
\begin{equation*}
m>s,1\leq p\leq \infty ,~~~1\leq q\leq \infty ,~~~y_{0}>0.
\end{equation*}%
The space $B_{p,q}^{s}(\Omega ;E)$ is $E$--valued Besov space, i.e., 
\begin{equation*}
\begin{array}{l}
B_{p,q}^{s}(\Omega ;E) \\ 
\vspace{-3mm} \\ 
\displaystyle=\left\{ 
\begin{array}{c}
f:f\in L_{p}(\Omega ;E),~\Vert f\Vert _{B_{p,q}^{s}(\Omega ;E)}~~ \\ 
=\sum\limits_{i=1}^{n}\left( \int\limits_{0}^{y_{0}}y^{-(sq+1)}\Vert \Delta
_{i}^{m}(y,\Omega )f(x)\Vert _{L_{p}(\Omega ;E)}^{q}dy\right) ^{\frac{1}{q}%
}<\infty%
\end{array}%
\right\} ,%
\end{array}%
\end{equation*}%
\begin{equation*}
\Vert f\Vert _{B_{p,\infty }^{s}(\Omega
;E)}~=~\sum\limits_{i=1}^{n}\sup\limits_{0<y<y_{0}}\frac{\Vert \Delta
_{i}^{m}(y,\Omega )f(x)\Vert _{L_{p}(\Omega ;E)}}{y^{s}},~~~1\leq p\leq
\infty ,~~~1\leq q<\infty .
\end{equation*}%
Let, 
\begin{equation*}
D^{\alpha }~=~D_{1}^{\alpha _{1}}D_{2}^{\alpha _{2}}\cdots D_{n}^{\alpha
_{n}},~~~D_{k}^{i}~=~\left( \frac{\partial }{\partial x_{k}}\right) ^{i}.
\end{equation*}%
and 
\begin{equation*}
\begin{array}{l}
\displaystyle B_{p,q}^{l,s}(\Omega ;E_{0},E)~=~\left\{ u:u\in
B_{p,q}^{s}(\Omega ;E_{0}),~~D_{k}^{l}u\in B_{p,q}^{s}(\Omega
;E),~~~k=1,2,...,n\right\} , \\ 
\vspace{-3mm} \\ 
\displaystyle\Vert u\Vert _{B_{p,q}^{l,s}(\Omega ;E_{0,}E)}~=~\Vert u\Vert
_{B_{p,q}^{s}(\Omega ;E_{0})}+\sum\limits_{k=1}^{n}\left\Vert
D_{k}^{l}u\right\Vert _{B_{p,q}^{s}(\Omega E)}~<~\infty .%
\end{array}%
\end{equation*}

$C^{(m)}(\Omega ;E)$ denotes the space of $E$--valued uniformly bounded and $%
m$--times continuously differentiable functions on $\Omega .$ \vspace{3mm}

$\left[ 7,\text{ \textbf{Lemma 2.3}}\right] .$ Let $\lambda \in S_{\varphi }$
and $\mu \in S_{\psi }$ where $\varphi +\psi <\pi \mathbf{.}$ Then there
exist $M>0$ such that $\left\vert \lambda +\mu \right\vert \geq M\left(
\left\vert \lambda \right\vert +\left\vert \mu \right\vert \right) .$

\vspace{3mm}

\section*{2. Youngs type Fourier multipliers}

Here, we shall study Fourier multiplier theorems from $B_{q_{1},r}^{s}$ to $%
B_{q_{2},r}^{s}$ for $\frac{1}{q_{2}}=\frac{1}{q_{1}}-\frac{1}{\eta ^{\prime
}}$ and $1<q_{1}<\eta ^{\prime }\leq \infty $ in the highlight of [10]. In
this section, $X$ and $Y$ are Banach spaces over the field $C$ and $X^{\ast
} $ is the dual space of $X.$ The space $B(X,Y)$ of bounded linear operators
from $X$ to $Y$ is endowed with the usual uniform operator topology. $N_{0}$
is the set of natural numbers containing zero.

All the basic properties of $F$ and $F^{-1}$ that hold in the scalar--valued
case also hold in vector--valued case; however, the Housdorff--Young
inequality need not hold. Therefore, we need to define Banach spaces that
was introduced by Peetre [11].

\vspace{3mm}

\textbf{Definition 2.1.} Let $X$ be a Banach space and $1\leq p\leq 2.$ We
say $X$ has Fourier type $p$ if 
\begin{equation*}
\Vert Ff\Vert _{L_{p^{\prime }}(R^{N},X)}~\leq ~C\Vert f\Vert
_{L_{p}(R^{N},X)}~~~\mbox{for each}~~~f\in S(R^{N},X),
\end{equation*}%
where $\frac{1}{p}+\frac{1}{p^{\prime }}=1,~F_{p,N}(X)$ is the smallest $%
C\in \lbrack 0,\infty ].$

\vspace{3mm}

$\left[ 10,\text{ \textbf{Proposition} }2.3\right] $ Let $X$ be a Banach
space with Fourier type $p\in \lbrack 1,2]$ and $p\leq q\leq p^{\prime }.$
Then $X^{\ast }$ and $L_{q}(R^{N},X)$ also have Fourier type provided both
are with the same constant $F_{p,N}(X).$

\hbox{\vrule
height7pt width5pt}~\hspace{3mm}

We shall use Fourier analytic definition of Besov spaces in this section.
Therefore, we need to consider some subsets $\{J_{k}\}_{k=0}^{\infty }$ and $%
\{I_{k}\}_{k=0}^{\infty }$ of $R^{N}$, where 
\begin{equation*}
J_{0}~=~\left\{ t\in R^{N}:~|t|\leq 1\right\} ,~~~J_{k}~=~\left\{ t\in
R^{N}:~2^{k-1}\leq |t|\leq 2^{k}\right\} ~~~\mbox{for}~~~k\in N
\end{equation*}%
and 
\begin{equation*}
I_{0}~=~\left\{ t\in R^{N}:~|t|\leq 2\right\} ,~~~I_{k}~=~\left\{ t\in
R^{N}:~2^{k-1}\leq |t|\leq 2^{k+1}\right\} ~~~\mbox{for}~~~k\in N.
\end{equation*}%
Next, we define the unity $\{\varphi _{k}\}_{k\in N_{0}}$ of functions from $%
S(R^{N},R).$ Let $\psi \in S(R,R)$ be nonnegative function with support in $%
[2^{-1},2],$ which satisfies 
\begin{equation*}
\sum\limits_{k=-\infty }^{\infty }\psi (2^{-k}s)~=~1~~~\mbox{for}~~~s\in
R\backslash \{0\}
\end{equation*}%
and 
\begin{equation*}
\varphi _{k}(t)~=~\psi (2^{-k}|t|),~~~\varphi
_{0}(t)~=~1-\sum\limits_{k=1}^{\infty }\varphi _{k}(t)~~~\mbox{for}~~~t\in
R^{N}.
\end{equation*}%
Later, we will need the following useful properties: 
\begin{equation*}
\begin{array}{l}
\mbox{supp}~\varphi _{k}~\subset ~\bar{I}_{k}~~~\mbox{for each}~~~k\in N_{0},
\\ 
\vspace{-3mm} \\ 
\varphi _{k}~\equiv ~0~~~\mbox{for each}~~~k<0, \\ 
\vspace{-3mm} \\ 
\displaystyle\sum\limits_{k=0}^{\infty }\varphi _{k}(s)~~~\mbox{for each}%
~~~s\in R^{N}, \\ 
\vspace{-3mm} \\ 
J_{m}\cap \mbox{supp}~\varphi _{k}~=~\emptyset ~~~\mbox{if}~~|m-k|>1, \\ 
\vspace{-3mm} \\ 
\varphi _{k-1}(s)+\varphi _{k}(s)+\varphi _{k+1}(s)~=~1~~~\mbox{for each}%
~~~s\in \mbox{supp}~\varphi _{k},~~~k\in N_{0}.%
\end{array}%
\end{equation*}

\vspace{3mm}

\textbf{Definition 2.3.} Let $1\leq q\leq r\leq \infty $ and $s\in R.$ The
Besov space is the set of all functions $f\in S^{\prime }(R^{N},X)$ for
which 
\begin{equation*}
\begin{array}{lll}
\Vert f\Vert _{B_{q,r}^{s}(R^{N},X)}: & = & \displaystyle\left\Vert
2^{ks}\left\{ (\check{\varphi}_{k}\ast f)\right\} _{k=0}^{\infty
}\right\Vert _{l_{r}(L_{q}(R^{N},X))} \\ 
\vspace{-3mm} &  &  \\ 
& \equiv & \displaystyle\left\{ 
\begin{array}{ll}
\displaystyle\left[ \sum\limits_{k=0}^{\infty }2^{ksr}\Vert \check{\varphi}%
_{k}\ast f\Vert _{L_{q}(R^{N},X)}^{r}\right] ^{\frac{1}{r}} & \mbox{if}%
~~r\neq \infty \\ 
\vspace{-3mm} &  \\ 
\displaystyle\sup\limits_{k\in N_{0}}\left[ 2^{ks}\Vert \check{\varphi}%
_{k}\ast f\Vert _{L_{q}(R^{N},X)}\right] & \mbox{if}~~r=\infty%
\end{array}%
\right.%
\end{array}%
\end{equation*}%
is finite; here $q$ and $s$ are main and smoothness indexes respectively.

To prove multiplier theorems, we will later need following results.

\vspace{3mm}

\textbf{\ }$\left[ 10,\text{ \textbf{Theorem 3.1}}\right] $ Let $X$ be a
Banach space with the Fourier type $p\in \lbrack 1,2].$ Let $1\leq
q<p^{\prime },~1\leq r\leq \infty $ and $s\geq $ $N\left( \frac{1}{q}-\frac{1%
}{p^{\prime }}\right) .$ Then, there exists a constant $C$ depending only on 
$F_{p,N}(X), $ so that if $f\in B_{p,r}^{s}(R^{N},X),$ 
\begin{equation*}
\left\Vert \left\{ \hat{f}.\chi _{J_{m}}\right\} _{k=0}^{\infty }\right\Vert
_{l_{r}(L_{q}(R^{N},X))}~\leq ~C\Vert f\parallel _{B_{p,r}^{s}(R^{N},X)}.
\end{equation*}

\vspace{3mm}

Note that $\left[ 10,\text{ \textbf{Theorem 3.1}}\right] $ remains valid if
Fourier transform is replaced by the inverse Fourier transform. Choosing $%
r=1 $ and $s=N\left( \frac{1}{q}-\frac{1}{p^{\prime }}\right) $ we obtain
the following corollary.

\vspace{3mm}

\textbf{Corollary 2.5.} Let $X$ be a Banach space having Fourier type $p\in
\lbrack 1,2]$ and $1\leq q<p^{\prime }.$ Then the Fourier transform defines
bounded operator 
\begin{equation*}
F:~B_{p,1}^{N\left( \frac{1}{q}-\frac{1}{p^{\prime }}\right)
}(R^{N},X)~\rightarrow ~L_{q}(R^{N},X).\eqno(1)
\end{equation*}

\vspace{3mm}

For a bounded measurable function $m:R^{N}\rightarrow B(X,Y),$ its
corresponding Fourier multiplier operator $T_{m}$ is defined as follows 
\begin{equation*}
T_{m}(f)~=~F^{-1}[m(\cdot )(Ff)(\cdot )].
\end{equation*}

In this section, we identify conditions on $m$, extending those of [10],
that 
\begin{equation*}
\Vert T_{m}f\Vert _{B_{q_{2},r}^{s}}~\leq ~C\Vert f\Vert
_{B_{q_{1},r}^{s}}~~~\mbox{for each}~~~f\in S(X).
\end{equation*}

\vspace{3mm}

\textbf{Definition 2.6.} Let $\left( E_{1}(R^{N},X)\text{,}%
E_{2}(R^{N},Y)\right) $ be one of the following systems, where $1\leq
q_{1},q_{2},~r\leq \infty $ and $s\in R$ 
\begin{equation*}
(L_{q_{1}}(X),L_{q_{2}}(Y))~~~\mbox{or}%
~~~(B_{q_{1},r}^{s}(X),B_{q_{2},r}^{s}(Y)).
\end{equation*}%
A bounded measurable function $m:R^{N}\rightarrow B(X,Y)$ is called a
Fourier multiplier from $E_{1}(X)$ to $E_{2}(Y)$ if there is a bounded
linear operator 
\begin{equation*}
T_{m}:~E_{1}(X)\rightarrow E_{2}(Y)
\end{equation*}%
such that 
\begin{equation*}
T_{m}(f)~=~F^{-1}[m(\cdot )(Ff)(\cdot )]~~~\mbox{for each}~~~f\in S(X),\eqno%
(2)
\end{equation*}%
\begin{equation*}
T_{m}~~~\mbox{is}~~~\sigma (E_{1}(X),E_{1}^{\ast }(X^{\ast }))~~~\mbox{to}%
~~~\sigma (E_{2}(Y),E_{2}^{\ast }(Y^{\ast }))~~~\mbox{continuous.}\eqno(3)
\end{equation*}%
The uniquely determined operator $T_{m}$ is the Fourier multiplier operator
induced by $m.$

\vspace{3mm}

\textbf{Remark 2.7.} If $T_{m}\in B(E_{1}(X),E_{2}(Y))$ and $T_{m}^{\ast }$
maps $E_{2}^{\ast }(Y^{\ast })$ into $E_{1}^{\ast }(X^{\ast })$ then $T_{m}$
satisfies the continuity condition (3).

\vspace{3mm}

\textbf{Lemma 2.8.} Let $k\in L_{\eta }(R^{N},B(X,Y))$, $\frac{1}{q_{2}}=%
\frac{1}{q_{1}}-\frac{1}{\eta ^{\prime }}$ and $\frac{1}{\eta }+\frac{1}{%
\eta ^{\prime }}=1$ for $1<q_{1}<\eta ^{\prime }\leq \infty .$ Assume there
exist $C_{0}$ so that for all $x\in X$ and $y^{\ast }\in Y^{\ast }$ 
\begin{equation*}
\Vert kx\Vert _{L_{\eta }(Y)}~\leq ~C_{0}\Vert x\Vert _{X}\eqno(4)
\end{equation*}%
and $C_{1}$ so that 
\begin{equation*}
\Vert k^{\ast }y^{\ast }\Vert _{L_{\eta }(X^{\ast })}~\leq ~C_{1}\Vert
y^{\ast }\Vert _{Y^{\ast }}.\eqno(5)
\end{equation*}%
Then the convolution operator 
\begin{equation*}
K:~L_{q_{1}}(R^{N},X)~\rightarrow ~L_{q_{2}}(R^{N},Y)
\end{equation*}%
defined by 
\begin{equation*}
(Kf)(t)~=~\int\limits_{R^{N}}k(t-s)f(s)ds~~~\mbox{for}~~~t\in R^{N}
\end{equation*}%
satisfies $\Vert K\Vert _{L_{q_{1}}\rightarrow L_{q_{2}}}\leq C.$

\vspace{2mm}

\textbf{Proof}. Taking into account that $1\leq \eta $ and by using general
Minkowski-Jessen inequality and $\left( 4\right) $ we obtain 
\begin{equation*}
\begin{array}{lll}
\Vert (Kf)(t)\Vert _{L_{\eta }(Y)} & \leq & \displaystyle\left[
\int\limits_{R^{N}}\left( \int\limits_{R^{N}}\Vert k(t-s)f(s)\Vert
_{Y}ds\right) ^{\eta }\,dt\right] ^{\frac{1}{\eta }} \\ 
&  &  \\ 
\vspace{-3mm} & \leq & \displaystyle\int\limits_{R^{N}}\left(
\int\limits_{R^{N}}\Vert k(t-s)f(s)\Vert _{Y}^{\eta }dt\right) ^{\frac{1}{%
\eta }}\,ds \\ 
&  &  \\ 
& = & \displaystyle\int\limits_{R^{N}}\Vert kf(s)\Vert _{L_{\eta }(Y)}\,ds
\\ 
&  &  \\ 
& \leq & \displaystyle~C_{0}\int\limits_{R^{N}}\Vert f(s)\Vert _{X}ds~\leq
~C_{0}\Vert f\Vert _{L_{1}(X)}.%
\end{array}%
\end{equation*}%
for all $f\in L_{1}(X).$ Thus, $\left\Vert K\right\Vert _{L_{1\rightarrow
}L_{\eta }}\leq C_{0}.$ Now, assume $f\in L_{\eta ^{\prime }}(X),~y^{\ast
}\in Y^{\ast }$ and $t\in R^{N}.$ Then, by using H\"{o}lder inequality and $%
\left( 5\right) $ we get 
\begin{equation*}
\begin{array}{lll}
|<y^{\ast },(Kf)(t)>_{Y}| & \leq & \displaystyle\int%
\limits_{R^{N}}|<k(t-s)^{\ast }y^{\ast },f(s)>_{X}|ds \\ 
\vspace{-3mm} &  &  \\ 
& \leq & \displaystyle\int\limits_{R^{N}}\left\vert k(t-s)^{\ast }y^{\ast
}f(s)\right\vert ds \\ 
\vspace{-3mm} &  &  \\ 
& \leq & \displaystyle\Vert k^{\ast }y^{\ast }\Vert _{L_{\eta }(X^{\ast
})}\Vert f(s)\Vert _{L_{\eta ^{\prime }}(X)} \\ 
&  &  \\ 
& \leq & \displaystyle C_{1}\Vert y^{\ast }\Vert \Vert f\Vert _{L_{\eta
^{\prime }}(X)}.%
\end{array}%
\end{equation*}

Hence, $\left\Vert K\right\Vert _{L_{\eta ^{\prime }}\rightarrow L_{\infty
}}\leq C_{1}.$ Now, Riesz-Thorin theorem implies that%
\begin{equation*}
\left\Vert K\right\Vert _{_{L_{q_{1}}\rightarrow L_{q_{2}}}}\leq C\text{ }
\end{equation*}%
for 
\begin{equation*}
\frac{1}{q_{1}}=1-\theta +\frac{\theta }{\eta ^{\prime }}\text{ and }\frac{1%
}{q_{2}}=\frac{1-\theta }{\eta }\text{ }
\end{equation*}%
where $0<\theta <1.$ Solving these equation we obtain 
\begin{equation*}
\left\Vert K\right\Vert _{_{L_{q_{1}}\rightarrow L_{q_{2}}}}\leq C
\end{equation*}%
for 
\begin{equation*}
\frac{1}{q_{2}}=\frac{1}{q_{1}}-\frac{1}{\eta ^{\prime }}.\text{ }
\end{equation*}%
\hbox{\vrule height7pt width5pt}

\vspace{3mm}

\textbf{Remark 2.9.} Let us note that if we choose $\eta ^{\prime }=\infty $
in Lemma 2.8, we obtain Lemma 4.5 of [10].

\vspace{3mm}

\textbf{Theorem 2.10.} Let $X$ and $Y$ be Banach spaces having Fourier type $%
p\in \lbrack 1,2]$. Then there is a constant $C$ depending only on $%
F_{p,N}(X)$ and $F_{p,N}(Y)$ so that if 
\begin{equation*}
m~\in ~B_{p,1}^{N\left( \frac{1}{\eta }-\frac{1}{p^{\prime }}\right)
}(R^{N},B(X,Y))
\end{equation*}%
then $m$ is a Fourier multiplier from $L_{q_{1}}(R^{N},X)$ to $%
L_{q_{2}}(R^{N},Y)$ with 
\begin{equation*}
\Vert T_{m}\Vert _{L_{q_{1}}(R^{N},X)\rightarrow L_{q_{2}}(R^{N},Y)}~\leq
~CM_{p,\eta }(m)~~~\mbox{for each}\eqno(6)
\end{equation*}%
where $\frac{1}{q_{2}}=\frac{1}{q_{1}}-\frac{1}{\eta ^{\prime }},$ $\frac{1}{%
\eta }+\frac{1}{\eta ^{\prime }}=1,$ $1<q_{1}<\eta ^{\prime }\leq \infty $
and 
\begin{equation*}
M_{p,\eta }(m)~=~\inf \left\{ \Vert m(a\cdot )\Vert _{B_{p,1}^{N\left( \frac{%
1}{\eta }-\frac{1}{p^{\prime }}\right) }(R^{N},B(X,Y))}:~a>0\right\}
\end{equation*}

\vspace{1mm}

\textbf{Proof.} The key point in this proof is the fact $\left( 1\right) $.
As in the proof of $[10,$Theorem $4.3]$ we assume in addition that $m\in
S\left( B\left( X,Y\right) \right) .$ Hence, $\check{m}\in S\left( B\left(
X,Y\right) \right) .$ Since, $F^{-1}\left[ m\left( a\cdot \right) x\right]
\left( s\right) =a^{-N}\check{m}\left( \frac{s}{a}\right) x$~choosing an
appropriate $a$ and using $\left( 1\right) $ we obtain 
\begin{equation*}
\begin{array}{lll}
\vspace{-3mm} &  & \displaystyle\left\Vert \check{m}\left( \cdot \right)
x\right\Vert _{L_{\eta }\left( Y\right) }=\left\Vert \left[ m\left( a\cdot
\right) x\right] ^{\vee }\right\Vert _{L_{\eta }\left( Y\right) } \\ 
&  &  \\ 
& \leq & \displaystyle C_{1}\left\Vert m\left( a\cdot \right) \right\Vert
_{B_{p,1}^{N\left( \frac{1}{\eta }-\frac{1}{p^{\prime }}\right) }}\left\Vert
x\right\Vert _{X} \\ 
&  &  \\ 
\vspace{-3mm} & \leq & \displaystyle2C_{1}M_{p,\eta }(m)\left\Vert
x\right\Vert _{X}%
\end{array}%
\end{equation*}%
~

where $C_{1}$ depends only on $F_{p,N}(Y).$~If $m\in S\left( B\left(
X,Y\right) \right) $ then $\left[ m(\cdot )^{\ast }\right] ^{\vee }=\left[ 
\check{m}(\cdot )\right] ^{\ast }\in S\left( B\left( Y^{\ast },X^{\ast
}\right) \right) $ and $M_{p,\eta }(m)=M_{p,\eta }(m^{\ast }).$ Thus, in a
similar manner as above, we have

$\ \ \ \ \ \ \ \ \ \ \ \ \ \ \ \ \ \ \ \ \ \ \ \ \ \ \ \ \ \ \ \ \ \ \ \ \ \
\ \ \ \ \ \ \ \ \ \ \ 
\begin{array}{lll}
\vspace{-3mm} &  & \displaystyle\left\Vert \left[ \check{m}\left( \cdot
\right) \right] ^{\ast }y^{\ast }\right\Vert _{L_{\eta }\left( Y\right)
}\leq 2C_{2}M_{p,\eta }(m)\left\Vert y^{\ast }\right\Vert _{Y^{\ast }}%
\end{array}%
$ for some constant $C_{2}$ depends on $F_{p,N}(X^{\ast }).$ Since, we have%
\begin{equation*}
\left\Vert \check{m}\left( \cdot \right) x\right\Vert _{L_{\eta }\left(
Y\right) }\leq 2C_{1}M_{p,\eta }(m)\left\Vert x\right\Vert _{X}
\end{equation*}%
and%
\begin{equation*}
\left\Vert \left[ \check{m}\left( \cdot \right) \right] ^{\ast }y^{\ast
}\right\Vert _{L_{\eta }\left( Y\right) }\leq 2C_{2}M_{p,\eta }(m)\left\Vert
y^{\ast }\right\Vert _{Y^{\ast }}
\end{equation*}%
by Lemma 2.8 we can conclude 
\begin{equation*}
\left( T_{m}f\right) \left( t\right) )=\dint\limits_{R^{N}}\check{m}\left(
t-s\right) f\left( s\right) ds
\end{equation*}%
satisfies%
\begin{equation*}
\left\Vert T_{m}f\right\Vert _{L_{q_{2}}\left( R^{N},Y\right) }\leq
CM_{p,\eta }(m)\left\Vert f\right\Vert _{L_{q_{1}}\left( R^{N},X\right) }
\end{equation*}%
for $\frac{1}{q_{2}}=\frac{1}{q_{1}}-\frac{1}{\eta ^{\prime }}$ where $%
1<q_{1}<\eta ^{\prime }\leq \infty .$ Now, taking into account the fact that 
$S\left( B\left( X,Y\right) \right) $ continuously embedded to $%
B_{p,1}^{N\left( \frac{1}{\eta }-\frac{1}{p^{\prime }}\right) }(B(X,Y))$ and
using the same reasoning as in the proof of $[10,$Theorem $4.3]$ one can
easily prove for the general case $m\in B_{p,1}^{N\left( \frac{1}{\eta }-%
\frac{1}{p^{\prime }}\right) }$ and that $T_{m}$ satisfies $\left( 2\right)
, $ $\left( 3\right) .$ \ \hbox{\vrule height7pt width5pt}

\vspace{3mm}

Let us note that choosing $\eta =1$ in the Theorem 2.10 , we obtain $[10,$%
Theorem $4.3]$. The next theorem is the extension of $[10,$Theorem $4.8].$

\vspace{3mm}

\textbf{Theorem 2.11.} Let $X$ and $Y$ be Banach spaces having Fourier type $%
p\in \lbrack 1,2]$ and $\frac{1}{q_{2}}=\frac{1}{q_{1}}-\frac{1}{\eta
^{\prime }},$ $\frac{1}{\eta }+\frac{1}{\eta ^{\prime }}=1$ for $%
1\,<q_{1}<\eta ^{\prime }\leq \infty .$ Then, there exist a constant $C$
depending only on $F_{p,N}(X)$ and $F_{p,N}(Y)$ so that if $%
m:R^{N}\rightarrow B(X,Y)$ satisfy 
\begin{equation*}
\varphi _{k}\cdot m\in B_{p,1}^{N\left( \frac{1}{\eta }-\frac{1}{p^{\prime }}%
\right) }(R^{N},B(X,Y))~~~\mbox{and}~~~M_{p,\eta }(\varphi _{k}\cdot m)~\leq
~A\eqno(7)
\end{equation*}%
then $m$ is Fourier multiplier from $B_{q_{1},r}^{s}(R^{N},X)$ to $%
B_{q_{2},r}^{s}(R^{N},Y)$ and $\Vert T_{m}\Vert \leq CA$ for each $s\in R$
and $r\in \lbrack 1,\infty ].$

\textbf{Proof. }Since $\varphi _{k}\cdot m\in B_{p,1}^{N\left( \frac{1}{\eta 
}-\frac{1}{p^{\prime }}\right) }(R^{N},B(X,Y)),$ Theorem 2.10 ensures that%
\begin{equation*}
\left\Vert T_{m\varphi _{k}}f\right\Vert _{L_{q_{2}}(R^{N},Y)}\leq
C~M_{p,\eta }(\varphi _{k}\cdot m)~\left\Vert f\right\Vert
_{L_{q_{1}}(R^{N},X)}\leq ~CA\left\Vert f\right\Vert _{L_{q_{1}}(R^{N},X)}.
\end{equation*}%
In the introduction we defined function $\psi _{k}=\varphi _{k-1}+\varphi
_{k}+\varphi _{k+1}$ that is equal to $1$ on supp$\varphi _{k}.$ Thus, 
\begin{equation*}
\begin{array}{lll}
\vspace{-3mm} &  & \displaystyle\left\Vert T_{m\psi _{k}}f\right\Vert
_{L_{q_{2}}\left( Y\right) }\leq \left\Vert T_{m\varphi _{k-1}}f\right\Vert
_{L_{q_{2}}\left( Y\right) } \\ 
&  &  \\ 
& + & \displaystyle\left\Vert T_{m\varphi _{k}}f\right\Vert
_{L_{q_{2}}\left( Y\right) }+\left\Vert T_{m\varphi _{k+1}}f\right\Vert
_{L_{q_{2}}\left( Y\right) } \\ 
&  &  \\ 
\vspace{-3mm} & \leq & \displaystyle3CA\left\Vert f\right\Vert
_{L_{q_{1}}(R^{N},X)}.%
\end{array}%
\end{equation*}%
Let $T_{0}:S(X)\rightarrow S^{\prime }(Y)$ be defined as%
\begin{equation*}
T_{0}f=F^{-1}\left[ m(\cdot )\left( Ff\right) (\cdot )\right] .
\end{equation*}%
From the proof of $[10,$ Theorem $4.3]$ we know that%
\begin{equation*}
\check{\varphi}_{k}\ast T_{0}f=T_{m\psi _{k}}\left( \check{\varphi}_{k}\ast
f\right) .
\end{equation*}%
Hence,%
\begin{equation*}
\begin{array}{lll}
\vspace{-3mm} &  & \displaystyle\left\Vert \check{\varphi}_{k}\ast
T_{0}f\right\Vert _{L_{q_{2}}\left( Y\right) }=\left\Vert T_{m\psi
_{k}}\left( \check{\varphi}_{k}\ast f\right) \right\Vert _{L_{q_{2}}\left(
Y\right) } \\ 
&  &  \\ 
& \leq & \displaystyle3CA\left\Vert f\right\Vert _{L_{q_{1}}(R^{N},X)}%
\end{array}%
\end{equation*}%
and%
\begin{equation*}
\begin{array}{lll}
\vspace{-3mm} &  & \displaystyle\dsum\limits_{k=0}^{\infty
}2^{ksr}\left\Vert \check{\varphi}_{k}\ast T_{0}f\right\Vert
_{L_{q_{2}}\left( Y\right) }^{r}\leq 3CA\dsum\limits_{k=0}^{\infty
}2^{ksr}\left\Vert f\right\Vert _{L_{q_{1}}(R^{N},X)}^{r}.%
\end{array}%
\end{equation*}%
Thus, we obtain that 
\begin{equation*}
\left\Vert T_{0}f\right\Vert _{B_{q_{2},r}^{s}(R^{N},Y)}\leq 3CA\left\Vert
f\right\Vert _{B_{q_{1},r}^{s}(R^{N},X)}
\end{equation*}%
for $\frac{1}{q_{2}}=\frac{1}{q_{1}}-\frac{1}{\eta ^{\prime }}$ where$%
1<q_{1}<\eta ^{\prime }\leq \infty .$ If $q,r<\infty $ then $\mathring{B}%
_{q,r}^{s}=B_{q,r}^{s}.$ Therefore, it remains only to show the weak
continuity condition $\left( 3\right) .$ Proof of the case $r=\infty $ and
the weak continuity condition $\left( 3\right) $ follows trivially from the
proof of $[10,$ Theorem $4.3].$ \ \hbox{\vrule height7pt width5pt}

\vspace{3mm}

To establish maximal regularity of DOEs and abstract embeddings we shall use
Theorem 2.11 in the next sections. However, it is not easy to check
assumptions of the Theorem 2.11 for multiplier functions. Therefore, we
prove a lemma that makes Theorem 2.11 more applicable.

\vspace{3mm}

\textbf{Lemma 2.12.} Let $N\left( \frac{1}{\eta }-\frac{1}{p^{\prime }}%
\right) <l\in N$ and $u\in \lbrack p,\infty ]$. If $m\in C^{l}(R^{N},B(X,Y))$
with 
\begin{eqnarray*}
\Vert D^{\alpha }m\left\vert _{I_{0}}\Vert _{L_{u}(B(X,Y))}\right. ~ &\leq
&~A,~~~\Vert D^{\alpha }m_{k}\left\vert _{I_{1}}\Vert
_{L_{u}(B(X,Y))}\right. ~\leq ~A, \\
m_{k}(\cdot ) &=&m(2^{k-1}\cdot ),
\end{eqnarray*}%
for each $\alpha \in N_{0}^{N},~|\alpha |\leq l$ and $k\in N$ then $m$
satisfies conditions of Theorem 2.11.

\textbf{Proof.} Using the fact $W_{p}^{l}(R^{N},B(X,Y))\subset
B_{p,1}^{N\left( \frac{1}{\eta }-\frac{1}{p^{\prime }}\right)
}(R^{N},B(X,Y)) $ for $N\left( \frac{1}{\eta }-\frac{1}{p^{\prime }}\right)
<l$ (see $\left[ 2\right] )$ one can prove this lemma in a similar manner as 
$[10,$ Lemma 4.10$]$.~~~~\hspace{3mm}\hbox{\vrule
height7pt width5pt}

\vspace{3mm}

The following corollary shows that classical Mikhlin condition implies
assumptions $(7)$ of Theorem 2.11.

\vspace{3mm}

\textbf{Corollary 2.13.} Let $X$ and $Y$ be Banach spaces having Fourier
type $p\in \lbrack 1,2]$ and $\frac{1}{q_{2}}=\frac{1}{q_{1}}-\frac{1}{\eta
^{\prime }},$ $\frac{1}{\eta }+\frac{1}{\eta ^{\prime }}=1$ for $%
1\,<q_{1}<\eta ^{\prime }\leq \infty .$ If $m\in C^{l}(R^{N},B(X,Y))$
satisfies 
\begin{equation*}
\left\Vert (1+|t|)^{|\alpha |}D^{\alpha }m(t)\right\Vert _{L_{\infty
}(R^{N},B(X,Y))}~\leq ~A\eqno(8)
\end{equation*}%
for each multi--index $\alpha $ with $|\alpha |\leq l=\left\lceil N\left( 
\frac{1}{\eta }-\frac{1}{p^{\prime }}\right) \right\rceil +1$ then $m$ is
Fourier multiplier from $B_{q_{1},r}^{s}(R^{N},X)$ to $%
B_{q_{2},r}^{s}(R^{N},Y)$ for each $s\in R$ and $r\in \lbrack 1,\infty ].$

\vspace{3mm}

\textbf{Proof.} The proof follows from $[10,$ Corollary 4.11$].$ Actually,
choosing $u=\infty $ in the Lemma 2.12 we get assertions of Corollary. \ 
\hbox{\vrule
height7pt width5pt}

\vspace{3mm}

\textbf{Remark 2.14. }It is known that, any Banach space has Fourier type $%
1. $ Therefore, in the Corollary 2.13 if $X$ and $Y$ are arbitrary Banach
spaces then $l=\left\lceil \frac{N}{\eta }\right\rceil +1$.

The next corollary shows that classical H\"{o}rmander condition implies
assumptions $(7)$ of Theorem 2.11.

\vspace{3mm}

\textbf{Corollary 2.15.} Assume $\frac{1}{q_{2}}=\frac{1}{q_{1}}-\frac{1}{%
\eta ^{\prime }}$ for $\frac{1}{\eta }+\frac{1}{\eta ^{\prime }}=1$ and $%
1\,<q_{1}<\eta ^{\prime }\leq \infty .$ Suppose $X$ and $Y$ have Fourier
type $p\in \lbrack 1,2]$ and $l=\left\lceil N\left( \frac{1}{\eta }-\frac{1}{%
p^{\prime }}\right) \right\rceil +1.$ If $m\in C^{l}(R^{N},B(X,Y))$
satisfies 
\begin{equation*}
\left[ \dint\limits_{\left\vert t\right\vert \leq 2}\left\Vert D^{\alpha
}m(t)\right\Vert ^{p}\right] ^{\frac{1}{p}}\leq ~A
\end{equation*}%
and 
\begin{equation*}
\left[ R^{-N}\dint\limits_{R\leq \left\vert t\right\vert \leq 4R}\left\Vert
D^{\alpha }m(t)\right\Vert ^{p}\right] ^{\frac{1}{p}}\leq ~AR^{-\left\vert
\alpha \right\vert }
\end{equation*}%
for each multi--index $\alpha $ with $|\alpha |\leq l$ then $m$ is Fourier
multiplier from $B_{q_{1},r}^{s}(R^{N},X)$ to $B_{q_{2},r}^{s}(R^{N},Y)$ for
each $s\in R$ and $r\in \lbrack 1,\infty ].$

\textbf{Proof. }Choosing $u=p$ in the Lemma 2.12 we get assertions of
corollary. \ \hbox{\vrule
height7pt width5pt}

\section*{3. Abstract Embeddings}

In the present section, using Corollary 2.13 we shall prove continuity of
the following embedding%
\begin{equation*}
D^{\alpha }:B_{q_{1},r}^{l,s}\left( R^{N};E\left( A\right) ,E\right) \subset
B_{q_{2},r}^{s}\left( R^{N};E\right) .
\end{equation*}

These type of embeddings very often used to establish maximal regularity for
differential operator equations see e.g $\left[ 14\right] ,$ $\left[ 15%
\right] $.

\vspace{3mm}

\textbf{Proposition 3.1}. Let $A$ be a positive operator in Banach space\ $E$%
, $\alpha =\left( \alpha _{1},\alpha _{2},...,\alpha _{N}\right) $ and $x=%
\frac{\left\vert \alpha \right\vert }{l}\leq 1$. Then, operator-function \ 

\begin{equation*}
\Psi \left( \xi \right) =\left( i\xi \right) ^{\alpha }A^{1-x}\left[
A+\theta \left( \xi \right) \right] ^{-1}
\end{equation*}%
is uniformly bounded and satisfies 
\begin{equation*}
\ \ \ \left\Vert \Psi \left( \xi \right) \right\Vert _{B\left( E\right)
}\leq C\eqno(10)
\end{equation*}%
for all $\xi \in R^{N},$ where 
\begin{equation*}
\theta =\theta \left( \xi \right) =\sum\limits_{k=1}^{N}\left\vert \xi
_{k}\right\vert ^{l}\in S\left( \varphi \right) .
\end{equation*}
\vspace{1mm} \ \ \ \ \ \ \ \textbf{Proof.}\ Since $\theta \left( \xi \right)
\in S\left( \varphi \right) ,$ for all $\varphi \in \left( 0\right. ,\left.
\pi \right] $ and $A$ is $\varphi $-positive in $E,$\ the operator $A+\theta
\left( \xi \right) $ is invertible$.$\ Let \ \ 
\begin{equation*}
u=\left[ A+\theta \left( \xi \right) \right] ^{-1}f.\eqno(11)
\end{equation*}%
Then 
\begin{equation*}
\left\Vert \Psi \left( \xi \right) f\right\Vert _{E}\leq \left\Vert
A^{1-x}u\right\Vert _{E}\left\vert \xi _{1}\right\vert ^{\alpha
_{1}}...\left\vert \xi _{N}\right\vert ^{\alpha N}.
\end{equation*}%
Using the Moment inequality for powers of positive operators, we get a
constant $C$ such that 
\begin{eqnarray*}
&&\left\Vert \Psi \left( \xi \right) f\right\Vert _{E} \\
&\leq &C\left\Vert Au\right\Vert ^{1-x}\left\Vert u\right\Vert
^{x}\left\vert \xi _{1}\right\vert ^{\alpha _{1}}...\left\vert \xi
_{N}\right\vert ^{\alpha _{N}}.
\end{eqnarray*}%
Then, applying Young inequality, which states that $ab\leq \frac{a^{k_{1}}}{%
k_{1}}+\frac{b^{k_{2}}}{k_{2}}$ for any positive real numbers\ $a,b$ and\ $%
k_{1},k_{2}$ with \ $\frac{1}{k_{1}}+\frac{1}{k_{2}}=1,$ to the product

\begin{equation*}
\left\Vert Au\right\Vert ^{1-x}\left[ \left\Vert u\right\Vert ^{x}\left\vert
\xi _{1}\right\vert ^{\alpha _{1}}...\left\vert \xi _{n}\right\vert ^{\alpha
_{N}}\right]
\end{equation*}%
we obtain 
\begin{equation*}
\begin{array}{lll}
\left\Vert \Psi \left( \xi \right) f\right\Vert _{E} & \leq & C\left\{
\left( 1-x\right) \left\Vert Au\right\Vert \right. \\ 
&  &  \\ 
& + & \displaystyle\left. x\left[ \left\vert \xi _{1}\right\vert \right] ^{%
\frac{\alpha _{1}}{x}}...\left\vert \xi _{N}\right\vert ^{\frac{\alpha _{N}}{%
x}}\ \ \left\Vert u\right\Vert \right\}%
\end{array}%
\eqno(12)
\end{equation*}

Since, 
\begin{equation*}
\ \sum\limits_{i=1}^{N}\frac{\alpha _{i}}{x}=l
\end{equation*}%
there exists a constant $M_{0}$ independent of $\xi ,$ such that

\begin{equation*}
\left\vert \xi _{1}\right\vert ^{\frac{\alpha _{1}}{x}}...\left\vert \xi
_{N}\right\vert ^{\frac{\alpha _{N}}{x}}\leq M_{0}\left(
1+\sum\limits_{k=1}^{N}\left\vert \xi _{k}\right\vert ^{l}\right)
\end{equation*}%
for all $\xi \in R^{N}.$ Combining above estimate with inequality $\left(
12\right) $ and using the fact that $x=\frac{\left\vert \alpha \right\vert }{%
l}\leq 1,$ we obtain

\begin{equation*}
\left\Vert \psi \left( \xi \right) f\right\Vert \leq C\left[ \left\Vert
Au\right\Vert +\sum\limits_{k=1}^{N}\left\vert \xi _{k}\right\vert
^{l}\left\Vert u\right\Vert +\left\Vert u\right\Vert \right] .
\end{equation*}

Then, with the help of $(11)$ we get

\begin{eqnarray*}
\left\Vert \psi \left( \xi \right) f\right\Vert &\leq &C\left[ \left\Vert A%
\left[ A+\theta \left( \xi \right) \right] ^{-1}f\right\Vert +\right. \\
&&\left. \sum\limits_{k=1}^{N}\left\vert \xi _{k}\right\vert ^{2l}\left\Vert 
\left[ A+\theta \left( \xi \right) \right] ^{-1}f\right\Vert +\left\Vert %
\left[ A+\theta \left( \xi \right) \right] ^{-1}f\right\Vert \right] .
\end{eqnarray*}%
Finally, using resolvent properties of positive operator $A$ we conclude

\begin{equation*}
\left\Vert \Psi \left( \xi \right) f\right\Vert _{E}\leq C_{0}\left\Vert
f\right\Vert _{E}
\end{equation*}%
for all $f\in E.$

\vspace{3mm}

\textbf{Proposition 3.2}. Let $A$ be a positive operator in Banach space\ $E$%
, $\alpha =\left( \alpha _{1},\alpha _{2},...,\alpha _{N}\right) $ and $x=%
\frac{\left\vert \alpha \right\vert +\sigma }{l}\leq 1$ for $\sigma
=\left\lceil N(1+\frac{1}{q_{2}}-\frac{1}{q_{1}})\right\rceil +1$. Then,
operator-function \ 

\begin{equation*}
\Psi \left( \xi \right) =\left\vert \xi \right\vert ^{\sigma }\left( i\xi
\right) ^{\alpha }A^{1-x}D^{\sigma }\left[ A+\theta \left( \xi \right) %
\right] ^{-1}
\end{equation*}%
is uniformly bounded.

\vspace{3mm}\textbf{Proof. }Since, $\frac{\left\vert \alpha \right\vert
+\sigma }{l}\leq 1$ we have 
\begin{eqnarray*}
\left\vert \xi \right\vert ^{\sigma }\left\vert i\xi \right\vert ^{\alpha }
&\leq &C\sum\limits_{j=1}^{N}|\xi _{j}|^{\sigma }\left\vert \xi
_{1}\right\vert ^{\alpha _{1}}\left\vert \xi _{2}\right\vert ^{\alpha
_{2}}\cdot \cdot \cdot \left\vert \xi _{N}\right\vert ^{\alpha _{N}} \\
&\leq &C_{o}\left( 1+\sum\limits_{k=1}^{N}\left\vert \xi _{k}\right\vert
^{l}\right) .
\end{eqnarray*}%
Thus, using above estimate and proof of Proposition 3.1 one can easily prove
assertion of this theorem.

\vspace{3mm}

\textbf{Theorem 3.3}. Let $E$ be a Banach space and $A$ \ be a $\varphi $%
-positive operator in $E,$ where $\varphi \in \left( 0\right. ,\left. \pi %
\right] .$ If $\alpha =\left( \alpha _{1},\alpha _{2},...,\alpha _{N}\right) 
$ and $x=\frac{\left\vert \alpha \right\vert +\sigma }{l}\leq 1$ for $\sigma
=\left\lceil N(1+\frac{1}{q_{2}}-\frac{1}{q_{1}})\right\rceil +1$ then the
following embedding 
\begin{equation*}
D^{\alpha }:B_{q_{1},r}^{l,s}\left( R^{N};E\left( A\right) ,E\right) \subset
B_{q_{2},r}^{s}\left( R^{N};E\left( A^{1-x}\right) \right)
\end{equation*}%
is continuous for$\frac{1}{q_{2}}=\frac{1}{q_{1}}-\frac{1}{\eta ^{\prime }},$
$\left( \frac{1}{\eta }+\frac{1}{\eta ^{\prime }}=1,1<q_{1}<\eta ^{\prime
}\leq \infty \right) $ and there exists a positive constant\ $C$, such that%
\begin{equation*}
\begin{array}{lll}
\left\Vert D^{\alpha }u\right\Vert _{B_{q_{2},r}^{s}\left( R^{N};E\left(
A^{1-x}\right) \right) } & \leq & C\left\Vert u\right\Vert
_{B_{q_{1},r}^{l,s}\left( R^{N};E\left( A\right) ,E\right) }%
\end{array}%
\eqno(13)
\end{equation*}%
for all $u\in B_{q_{1},r}^{l,s}\left( R^{N};E\left( A\right) ,E\right) .$

\ 

\textbf{Proof}. Since $A$ is constant and closed operator$,$ we have%
\begin{equation*}
\begin{array}{lll}
\left\Vert D^{\alpha }u\right\Vert _{B_{q_{2},r}^{s}\left( R^{N};E\left(
A^{1-x}\right) \right) } & = & \left\Vert A^{1-x}D^{\alpha }u\right\Vert
_{B_{q_{2},r}^{s}\left( R^{N};E\right) } \\ 
& \backsim & \left\Vert F^{-\shortmid }\left( i\xi \right) ^{\alpha
}A^{1-x}Fu\right\Vert _{B_{q_{2},r}^{s}\left( R^{N};E\right) }.%
\end{array}%
\eqno(14)
\end{equation*}
(The symbol $\backsim $ indicates norm equivalency). In a similar manner,
from definition of $B_{q_{1},r}^{l,s}\left( R^{N};E_{0},E\right) $ we have%
\begin{equation*}
\left\Vert u\right\Vert _{B_{q_{1},r}^{l,s}\left( R^{N};E_{0},E\right) }\sim
\left\Vert Au\right\Vert _{B_{q_{1},r}^{s}\left( R^{N};E\right)
}+\sum\limits_{k=1}^{N}\left\Vert F^{-1}\xi _{k}^{l}\hat{u}\right\Vert
_{B_{q_{1},r}^{s}\left( R^{N};E\right) }.\eqno(15)
\end{equation*}%
By virtue of above relations, it is sufficient to prove 
\begin{eqnarray*}
&&\left\Vert F^{-\shortmid }\left[ \left( i\xi \right) ^{\alpha }A^{1-x}\hat{%
u}\right] \right\Vert _{B_{q_{2},r}^{s}\left( R^{N};E\right) } \\
&\leq &C\left[ \left\Vert F^{-\shortmid }A\hat{u}\right\Vert
_{B_{q_{1},r}^{s}\left( R^{N};E\right) }+\sum\limits_{k=1}^{N}\left\Vert
F^{-\shortmid }\left( \xi _{k}^{l}\hat{u}\right) \right\Vert
_{B_{q_{1},r}^{s}\left( R^{N};E\right) }\right] .
\end{eqnarray*}%
\ Hence, the inequality (13) will be followed if we can prove the following
estimate \ 
\begin{equation*}
\left\Vert F^{-\shortmid }\left[ \left( i\xi \right) ^{\alpha }A^{1-x}\hat{u}%
\right] \right\Vert _{B_{q_{2},r}^{s}\left( R^{N};E\right) }\leq C\left\Vert
F^{-\shortmid }\left( \left[ A+I\theta \right] \hat{u}\right) \right\Vert
_{B_{q_{1},r}^{s}\left( R^{N};E\right) }\eqno(16)
\end{equation*}%
for all $u\in B_{q_{1},r}^{l,s}\left( R^{N};E\left( A\right) ,E\right) ,$
where 
\begin{equation*}
\theta =\theta \left( \xi \right) =\sum\limits_{k=1}^{N}\left\vert \xi
_{k}\right\vert ^{l}\in S\left( \varphi \right) .
\end{equation*}%
Let us express the left hand side of (16) as follows 
\begin{eqnarray*}
&&\left\Vert F^{-\shortmid }\left[ \left( i\xi \right) ^{\alpha }A^{1-x}\hat{%
u}\right] \right\Vert _{B_{q_{2},r}^{s}\left( R^{N};E\right) } \\
&=&\left\Vert F^{-\shortmid }\left( i\xi \right) ^{\alpha }A^{1-x}\left[
(A+I\theta \right] ^{-1}\left[ \left( A+I\theta \right) \right] \hat{u}%
\right\Vert _{B_{q_{2},r}^{s}\left( R^{N};E\right) }
\end{eqnarray*}%
(Since\ $A$ is the positive operator in\ $E$ and $\theta \left( \xi \right)
\in S\left( \varphi \right) ,$ $\left[ (A+I\theta \right] ^{-1}$ exists ).
From Corollary 2.13 we know that%
\begin{equation*}
\begin{array}{lll}
\left\Vert F^{-\shortmid }\left( i\xi \right) ^{\alpha }A^{1-x}\left[
(A+I\theta \right] ^{-1}\left[ \left( A+I\theta \right) \right] \hat{u}%
\right\Vert _{B_{q_{2},r}^{s}\left( R^{N};E\right) } & \leq &  \\ 
C\left\Vert F^{-\shortmid }\left( \left[ A+I\theta \right] \hat{u}\right)
\right\Vert _{B_{q_{1},r}^{s}\left( R^{N};E\right) } &  & 
\end{array}%
\eqno(17)
\end{equation*}%
holds if $\frac{1}{q_{2}}=\frac{1}{q_{1}}-\frac{1}{\eta ^{\prime }}$ and
operator-function $\left( i\xi \right) ^{\alpha }A^{1-x}\left[ (A+\theta %
\right] ^{-1}$ satisfies (8) for each multi--index $\beta ,$ $|\beta |\leq
\left\lceil \frac{N}{\eta }\right\rceil +1.$ It is clear that 
\begin{equation*}
\begin{array}{rll}
\left\Vert (1+|\xi |)^{|\beta |}D^{\beta }\Psi \left( \xi \right)
\right\Vert _{L_{\infty }(B(E))} & \leq & \displaystyle\sum\limits_{k=0}^{|%
\beta |}\left\Vert |\xi |^{k}D^{\beta }\Psi \left( \xi \right) \right\Vert
_{L_{\infty }(B(E))} \\ 
\vspace{-3mm} &  & 
\end{array}%
\eqno(18)
\end{equation*}%
Therefore, it is enough to show 
\begin{equation*}
\begin{array}{rll}
\left\Vert |\xi |^{k}D^{\beta }\Psi \left( \xi \right) \right\Vert
_{L_{\infty }(B(E))} & \leq & C \\ 
\vspace{-3mm} &  & 
\end{array}%
\end{equation*}%
for $k=0,1,\cdot \cdot \cdot |\beta |$ and $|\beta |\leq \left\lceil N(1+%
\frac{1}{q_{2}}-\frac{1}{q_{1}})\right\rceil +1\leq N+1.$ The case $k=0$ and 
$|\beta |=0$ follows from Proposition 3.1 and the principal case $k=|\beta
|=\left\lceil N(1+\frac{1}{q_{2}}-\frac{1}{q_{1}})\right\rceil +1$ follows
from Proposition 3.2. For the sake of simplicity of calculations, we shall
prove only for the cases $|\beta |=1,~k=0$ and $k=1.$ Taking derivative of
operator function $\Psi \left( \xi \right) $ with respect to $\xi _{i}$ we
get 
\begin{equation*}
\begin{array}{rll}
\left\Vert \frac{\partial }{\partial \xi _{i}}\Psi \left( \xi \right)
\right\Vert _{L_{\infty }(B(E))} & \leq & \left\Vert I_{1}\right\Vert
_{B(E)}+\left\Vert I_{2}\right\Vert _{B(E)}+\left\Vert I_{3}\right\Vert
_{B(E)}%
\end{array}%
\end{equation*}%
where%
\begin{equation*}
\begin{array}{rll}
I_{1} & = & \left\vert \xi _{1}\right\vert ^{\alpha _{1}}\cdot \cdot \cdot
\left\vert \alpha _{i}\xi _{i}\right\vert ^{\alpha _{i}-1}\cdots \left\vert
\xi _{N}\right\vert ^{\alpha _{N}}A^{1-x}\left[ (A+\theta \right] ^{-1} \\ 
\vspace{-3mm}I_{2} & = & \left\vert \xi _{1}\right\vert ^{\alpha
_{1}}\left\vert \xi _{2}\right\vert ^{\alpha _{2}}\cdots \left\vert \xi
_{N}\right\vert ^{\alpha _{N}}A^{1-x}\sum\limits_{\substack{ k=1  \\ k\neq i 
}} ^{N}\left( \xi _{k}\right) ^{l}\left[ (A+\theta \right] ^{-2} \\ 
&  &  \\ 
I_{3} & = & \left\vert \xi _{1}\right\vert ^{\alpha _{1}}\left\vert \xi
_{2}\right\vert ^{\alpha _{2}}\cdots \left\vert \xi _{N}\right\vert ^{\alpha
_{N}}A^{1-x}l\left( \xi _{i}\right) ^{l-1}\left[ (A+\theta \right] ^{-2}.%
\end{array}%
\end{equation*}%
Since, the first term $I_{1}$ is almost same with $\Psi \left( \xi \right) $
one can easily show that it is uniformly bounded. By virtue of the facts
that $\Psi \left( \xi \right) $ is uniformly bounded, $A$ is closed and
positive operator in\ $E$ and $\theta \left( \xi \right) \in S\left( \varphi
\right) $ we obtain%
\begin{equation*}
\begin{array}{rll}
\left\Vert \vspace{-3mm}I_{2}\right\Vert _{B(E)} & = & \left\Vert \left\vert
\xi _{1}\right\vert ^{\alpha _{1}}\left\vert \xi _{2}\right\vert ^{\alpha
_{2}}\cdots \left\vert \xi _{N}\right\vert ^{\alpha _{N}}A^{1-x}\left[
(A+\theta \right] ^{-1}\sum\limits_{\substack{ k=1  \\ k\neq i }} ^{N}\left(
\xi _{k}\right) ^{l}\left[ (A+\theta \right] ^{-1}\right\Vert \\ 
&  &  \\ 
& = & \left\Vert \Psi \left( \xi \right) \sum\limits_{\substack{ k=1  \\ %
k\neq i }} ^{N}\left( \xi _{k}\right) ^{l}\left[ (A+\theta \right]
^{-1}\right\Vert \leq C\sum\limits_{\substack{ k=1  \\ k\neq i }}
^{N}\left\vert \left( \xi _{k}\right) \right\vert ^{l}\left(
1+\sum\limits_{k=1}^{N}\left\vert \xi _{k}\right\vert ^{l}\right) ^{-1} \\ 
&  &  \\ 
& \leq & C.%
\end{array}%
\end{equation*}%
Taking into account 
\begin{equation*}
l\left\vert \left( \xi _{i}\right) \right\vert ^{l-1}\leq
C\sum\limits_{k=1}^{N}\left\vert \xi _{k}\right\vert ^{l}
\end{equation*}%
and using the same reasoning as above one can easily prove that $I_{3}$ is
uniformly bounded. Hence,%
\begin{equation*}
\left\Vert \frac{\partial }{\partial \xi _{i}}\Psi \left( \xi \right)
\right\Vert _{L_{\infty }(B(E))}\leq C.
\end{equation*}%
In order to show

\begin{equation*}
\left\Vert |\xi |\frac{\partial }{\partial \xi _{i}}\Psi \left( \xi \right)
\right\Vert _{L_{\infty }(B(E))}\leq C
\end{equation*}%
it is sufficient to prove $|\xi |I_{j}$ are bounded for $j=1,2,3.$ It is
clear that,%
\begin{equation*}
\begin{array}{rll}
\left\Vert |\xi |I_{1}\right\Vert _{B(E)} & \leq & \sum\limits_{j=1}^{N}%
\left\Vert |\xi _{j}|\left\vert \xi _{1}\right\vert ^{\alpha _{1}}\cdot
\cdot \cdot \left\vert \alpha _{i}\xi _{i}\right\vert ^{\alpha _{i}-1}\cdots
\left\vert \xi _{N}\right\vert ^{\alpha _{N}}A^{1-x}\left[ (A+\theta \right]
^{-1}\right\Vert _{B(E)}%
\end{array}%
\end{equation*}%
Since, there exist a constant $M_{0}$ independent of $\xi $ such that

\begin{equation*}
\begin{array}{rll}
|\xi _{j}|\left\vert \xi _{1}\right\vert ^{\alpha _{1}}\cdot \cdot \cdot
\left\vert \xi _{i}\right\vert ^{\alpha _{i}-1}\cdots \left\vert \xi
_{N}\right\vert ^{\alpha _{N}} & \leq & M_{0}\left(
1+\sum\limits_{k=1}^{N}\left\vert \xi _{k}\right\vert ^{l}\right)%
\end{array}%
\end{equation*}%
for all $j=0,1\cdot \cdot \cdot ,N,$ using same arguments as in Proposition
3.1 one can easily show that%
\begin{equation*}
\begin{array}{rll}
\left\Vert |\xi |I_{1}\right\Vert _{B(E)} & \leq & C.%
\end{array}%
\end{equation*}%
Similarly, by virtue of Proposition 3.1 and using the same techniques, one
can establish uniformly boundedness of $|\xi |I_{2}$ and $|\xi |I_{3}.$
Hence, operator functions $\frac{\partial }{\partial \xi _{i}}\Psi \left(
\xi \right) $ and $|\xi |\frac{\partial }{\partial \xi _{i}}\Psi \left( \xi
\right) $ are uniformly bounded multipliers. Other cases can be proved
analogously. \ \hbox{\vrule height7pt width5pt}

The next embedding theorem arises in the investigation of COE's.

\vspace{3mm}

\textbf{Theorem 3.5}. Let $E$ be a Banach space and $A$ \ be a $\varphi $%
-positive operator in $E,$ where $\varphi \in \left( 0\right. ,\left. \pi %
\right] .$ If $a_{\alpha }\in L_{1}(R)$, $\alpha =\left( \alpha _{1},\alpha
_{2},...,\alpha _{N}\right) $ and $x=\frac{\left\vert \alpha \right\vert
+\sigma }{l}\leq 1$ for $\sigma =\left\lceil N(1+\frac{1}{q_{2}}-\frac{1}{%
q_{1}})\right\rceil +1$ then the following embedding 
\begin{equation*}
a_{\alpha }\ast D^{\alpha }:B_{q_{1},r}^{l,s}\left( R^{N};E\left( A\right)
,E\right) \subset B_{q_{2},r}^{s}\left( R^{N};E\left( A^{1-x}\right) \right)
\end{equation*}%
is continuous for$\frac{1}{q_{2}}=\frac{1}{q_{1}}-\frac{1}{\eta ^{\prime }},$
$\left( \frac{1}{\eta }+\frac{1}{\eta ^{\prime }}=1,1<q_{1}<\eta ^{\prime
}\leq \infty \right) $ and there exists a positive constant\ $C$, such that 
\begin{equation*}
\begin{array}{lll}
\left\Vert a_{\alpha }\ast D^{\alpha }u\right\Vert _{B_{q_{2},r}^{s}\left(
R^{N};E\left( A^{1-x}\right) \right) } & \leq & C_{a_{\alpha }}\left\Vert
u\right\Vert _{B_{q_{1},r}^{l,s}\left( R^{N};E\left( A\right) ,E\right) }%
\end{array}%
\eqno(19)
\end{equation*}%
for all $u\in B_{q_{1},r}^{l,s}\left( R^{N};E\left( A\right) ,E\right) .$ 
\vspace{1mm}

\textbf{Proof}. Really, using Theorem 3.3 and Youngs inequality we get%
\begin{equation*}
\begin{array}{lll}
\left\Vert a_{\alpha }\ast D^{\alpha }u\right\Vert _{B_{q_{2},r}^{s}\left(
R^{N};E\left( A^{1-x}\right) \right) } & \leq & \displaystyle\left\Vert
a_{\alpha }\right\Vert _{L_{1}(R)}\left\Vert D^{\alpha }u\right\Vert
_{B_{q_{2},r}^{s}\left( R^{N};E\left( A^{1-x}\right) \right) } \\ 
\vspace{-3mm} &  &  \\ 
& \leq & \displaystyle C_{a_{\alpha }}\left\Vert u\right\Vert
_{B_{q_{1},r}^{l,s}\left( R^{N};E\left( A\right) ,E\right) }.%
\end{array}%
\end{equation*}%
\hbox{\vrule height7pt width5pt}

\vspace{3mm}

\textbf{Result 3.6}. Let all conditions of Theorem 3.3 hold and $\left\vert
\alpha \right\vert =l-\sigma $. Then for all $u\in B_{q_{1},s}^{l,s}\left(
R^{N};E\left( A\right) ,E\right) $ we have 
\begin{equation*}
D^{\alpha }:B_{q_{1},r}^{l,s}\left( R^{N};E\left( A\right) ,E\right) \subset
B_{q_{2},r}^{s}\left( R^{N};E\right) .
\end{equation*}%
Indeed, choosing $\left\vert \alpha \right\vert =l-\sigma $ in Theorem 3.3
we get $x=1$ that implies desired result.

\vspace{3mm}

\textbf{Result 3.7}. Let all conditions of Theorem 3.5 hold and $\left\vert
\alpha \right\vert =l-\sigma $. Then for all $u\in B_{q_{1},s}^{l,s}\left(
R^{N};E\left( A\right) ,E\right) $ we have 
\begin{equation*}
a_{\alpha }\ast D^{\alpha }:B_{q_{1},r}^{l,s}\left( R^{N};E\left( A\right)
,E\right) \subset B_{q_{2},r}^{s}\left( R^{N};E\right) .
\end{equation*}

\hbox{\vrule height7pt width5pt}

\section*{4. Differential--operator equations}

The main aim of this section is to establish $B_{q_{1},r}^{s}\rightarrow
B_{q_{2},r}^{s}$ regularity of the following elliptic DOE 
\begin{equation*}
(Q+\lambda )u~=~\sum\limits_{|\alpha |\leq 2l}a_{\alpha }D^{\alpha
}u+A_{\lambda }u~=~f,\eqno(20)
\end{equation*}%
where $A_{\lambda }=A+\lambda $ is a possible unbounded operator in $E.$

\vspace{3mm}

\textbf{Theorem 4.1.} Suppose $\frac{1}{q_{2}}=\frac{1}{q_{1}}-\frac{1}{\eta
^{\prime }}$ $\left( \text{ for }\frac{1}{\eta }+\frac{1}{\eta ^{\prime }}=1%
\text{ and }1<q_{1}<\eta ^{\prime }\leq \infty \right) $ and the following
conditions hold:

\vspace{2mm}

\noindent (1) $E$ is a Banach space;

\noindent (2) $A$ is a $\varphi $--positive operator in $E$ with $\varphi
\in \lbrack 0,\pi )$ and 
\begin{equation*}
\begin{array}{lll}
L(\xi ) & = & \displaystyle\sum\limits_{|\alpha |\leq 2l}a_{\alpha }(i\xi
_{1})^{\alpha _{1}}(i\xi _{2})^{\alpha _{2}}\cdots (i\xi _{N})^{\alpha
_{N}}\in S(\varphi _{1}), \\ 
\vspace{-3mm} &  &  \\ 
|L(\xi )| & \geq & \displaystyle K\sum\limits_{k=1}^{N}|\xi
_{k}|^{2l},~~~\xi \in R^{N},~~~\varphi _{1}+\varphi <\pi .%
\end{array}%
\end{equation*}%
Then for all $f\in B_{q_{1},r}^{s}(R^{N};E),~r\in \lbrack 1,\infty ]$, $%
\lambda \in S(\varphi )$ the equation $(20)$ has a unique solution $u(x)\in
B_{q_{2},r}^{2l,s}(R^{N};E(A),E)$ and coercive uniform estimate 
\begin{equation*}
\sum\limits_{|\alpha |\leq 2l}|\lambda |^{1-\frac{|\alpha |}{2l}}\left\Vert
D^{\alpha }u\right\Vert _{B_{q_{2},r}^{s}(R^{N};E)}+\Vert Au\Vert
_{B_{q_{2},r}^{s}(R^{N};E))}~\leq ~C\Vert f\Vert _{B_{q_{1},r}^{s}(R^{N};E)}%
\eqno(21)
\end{equation*}%
holds.

\vspace{1mm}

\textbf{Proof.} Applying Fourier transform to equation $(20)$, we obtain 
\begin{equation*}
\lbrack L(\xi )+A+\lambda ]\hat{u}(\xi )~=~\hat{f}(\xi ).\eqno(22)
\end{equation*}

Since $L(\xi )\in S(\varphi _{1})$ for all $\xi \in R^{N}$ and $A$ is a
positive operator, solutions are of the form 
\begin{equation*}
u(x)~=~F^{-1}[A+\lambda +L(\xi )]^{-1}\hat{f}.\eqno(23)
\end{equation*}

By using $(23)$, we have 
\begin{equation*}
\begin{array}{lll}
\Vert Au\Vert _{B_{q_{2},r}^{s}(R^{N};E)} & = & \displaystyle\left\Vert
F^{-1}A[A+(\lambda +L(\xi ))]^{-1}\hat{f}\right\Vert
_{B_{q_{2},r}^{s}(R^{N};E)} \\ 
\vspace{-3mm} &  &  \\ 
\Vert D^{\alpha }u\Vert _{B_{q_{2},r}^{s}(R^{N};E)} & = & \displaystyle%
\left\Vert F^{-1}[L_{1}\hat{f}]\right\Vert _{B_{q_{2},r}^{s}(R^{N};E)},%
\end{array}%
\end{equation*}%
where 
\begin{equation*}
L_{1}(\xi )~=~(i\xi _{1})^{\alpha _{1}}(i\xi _{2})^{\alpha _{2}}\cdots (i\xi
_{N})^{\alpha _{N}}[A+(\lambda +L(\xi ))]^{-1}.
\end{equation*}%
Hence, it suffices to show that operator--functions 
\begin{equation*}
\sigma _{1\lambda }(\xi )~=~A[A+(\lambda +L(\xi ))]^{-1},~~~\sigma
_{2\lambda }(\xi )~=~\sum\limits_{|\alpha |\leq 2l}|\lambda |^{1-\frac{%
|\alpha |}{2l}}L_{1}(\xi )
\end{equation*}%
are uniformly bounded multipliers from $B_{q_{1},r}^{s}(R^{N};E)$ to $%
B_{q_{2},r}^{s}(R^{N};E).$ In order to use Corollary 2.13, we have to show
that $\sigma _{j\lambda }\in C^{l}(R^{N},B(E))$ and there exists a constant $%
C>0$ such that 
\begin{equation*}
\left\Vert (1+|\xi |)^{|\beta |}D^{\beta }\sigma _{j\lambda }(\xi
)\right\Vert _{L_{\infty }(R^{n},B(E))}~\leq ~C
\end{equation*}%
for each multi--index $\beta $ with $|\beta |\leq \left\lceil \frac{N}{\eta }%
\right\rceil +1\leq N+1.$ For this, from the resolvent property of positive
operator $A$, we have 
\begin{equation*}
\begin{array}{lll}
\Vert \sigma _{1\lambda }(\xi )\Vert _{B(E)} & = & \displaystyle\left\Vert
A[A+(\lambda +L(\xi ))]^{-1}\right\Vert _{B(E)} \\ 
\vspace{-3mm} &  &  \\ 
& = & \displaystyle\left\Vert I-(\lambda +L(\xi ))[A+(\lambda +L(\xi
))]^{-1}\right\Vert _{B(E)} \\ 
\vspace{-3mm} &  &  \\ 
& \leq & \displaystyle1+|\lambda +L(\xi )|\left\Vert [A+(\lambda +L(\xi
))]^{-1}\right\Vert _{B(E)} \\ 
\vspace{-3mm} &  &  \\ 
& \leq & 1+M\Vert \lambda +L(\xi )|(1+|\lambda +L(\xi )|)^{-1}~\leq ~1+M.%
\end{array}%
\eqno(24)
\end{equation*}%
Now let us consider $\sigma _{2\lambda }.$ It is clear that%
\begin{equation*}
\Vert \sigma _{2\lambda }(\xi )f\Vert _{E}~=~\left\Vert |\lambda |^{1-\frac{%
|\alpha |}{2l}}L_{1}(\xi )f\right\Vert _{E}~\leq ~|\lambda |\left( |\xi
_{1}||\lambda |^{\frac{-1}{2l}}\right) ^{\alpha _{1}}\cdots \left( |\xi
_{N}||\lambda |^{\frac{-1}{2l}}\right) ^{\alpha _{N}}\left\Vert [A+(\lambda
+L(\xi ))]^{-1}f\right\Vert _{E}.\eqno(25)
\end{equation*}%
Then, by using the well--known inequality 
\begin{equation*}
|\xi _{1}|^{\alpha _{1}}|\xi _{2}|^{\alpha _{2}}\cdots |\xi _{N}|^{\alpha
_{N}}~\leq ~C\left( 1+\sum\limits_{k=1}^{N}|\xi _{k}|^{2l}\right)
,~~~|\alpha |\leq 2l\eqno(26)
\end{equation*}%
[7, Lemma 2.3], (25), the positivity of operator $A$ and ellipticity of $L,$
we obtain that 
\begin{equation*}
\begin{array}{lll}
\Vert \sigma _{2\lambda }(\xi )f\Vert _{E} & \leq & \displaystyle CM\left(
|\lambda |+\sum\limits_{k=1}^{N}|\xi _{k}|^{2l}\right) (1+|\lambda +L(\xi
)|)^{-1}\left\Vert f\right\Vert _{E} \\ 
\vspace{-3mm} &  &  \\ 
& \leq & C_{0}\left( |\lambda |+\sum\limits_{k=1}^{N}|\xi _{k}|^{2l}\right)
(1+|\lambda |+|L(\xi )|)^{-1}\Vert f\Vert _{E}~\leq ~C_{1}\Vert f\Vert _{E}%
\end{array}%
\eqno(27)
\end{equation*}%
for all $\lambda \in S(\varphi _{1})$ and $\xi \in R^{N}$.

It is clear that 
\begin{equation*}
\begin{array}{rll}
\left\Vert (1+|\xi |)^{|\beta |}D^{\beta }\sigma _{j\lambda }(\xi
)\right\Vert _{L_{\infty }(B(E)))} & \leq & \displaystyle\sum%
\limits_{k=0}^{|\beta |}\left\Vert |\xi |^{k}D^{\beta }\sigma _{j\lambda
}(\xi )\right\Vert _{L_{\infty }(B(E)))} \\ 
\vspace{-3mm} &  & 
\end{array}%
\eqno(28)
\end{equation*}%
Without any loss of generality, we shall prove (28) for $\sigma _{1\lambda
}(\xi ).$ For the sake of simplicity of calculations, we shall prove for
cases $|\beta |=1,~k=0$ and $k=1.$ Let us first recall the following
well--known inequality 
\begin{equation*}
\prod\limits_{k=1}^{N}|\xi _{k}|^{\alpha _{k}}~\leq ~M\left[
1+\sum\limits_{k=1}^{N}|\xi _{k}|^{2l}\right] ~~~\mbox{for}~~~|\alpha |\leq
2l.\eqno(29)
\end{equation*}%
Since 
\begin{equation*}
\frac{\partial }{\partial \xi _{k}}\sigma _{1\lambda }(\xi )~=~-A[A+(\lambda
+L(\xi ))]^{-2}\sum\limits_{|\alpha |\leq 2l}a_{\alpha }(i\xi _{1})^{\alpha
_{1}}(i\xi _{2})^{\alpha _{2}}\cdots (i\xi _{k})^{\alpha _{k}-1}\cdots (i\xi
_{N})^{\alpha _{N}}
\end{equation*}%
by using (29), we obtain 
\begin{equation*}
\begin{array}{lll}
\left\Vert \frac{\partial }{\partial \xi _{k}}\sigma _{1\lambda }(\xi
)\right\Vert _{B(E)} & \leq & 
\begin{array}{l}
\displaystyle C\sum\limits_{|\alpha |\leq 2l}|\xi _{1}|^{\alpha _{1}}|\xi
_{2}|^{\alpha _{2}}\cdots |\xi _{k}|^{\alpha _{k}-1} \\ 
\cdots |\xi _{N}|^{\alpha _{N}}\left\Vert A[A+(\lambda +L(\xi
))]^{-2}\right\Vert _{B(E)}%
\end{array}
\\ 
\vspace{-3mm} &  &  \\ 
& \leq & \displaystyle C\left[ 1+\sum\limits_{k=1}^{N}|\xi _{k}|^{2l}\right]
\left\Vert A[A+(\lambda +L(\xi ))]^{-2}\right\Vert _{B(E)}%
\end{array}%
\end{equation*}%
where $|\alpha |-1\leq 2l.$ Then by using [7, Lemma 2.3], resolvent
properties of positive operator $A$ and ellipticity of $L,$ for all $\lambda
\in S(\varphi _{1}),~\xi \in R^{N}$, we have 
\begin{equation*}
\begin{array}{lll}
\left\Vert \frac{\partial }{\partial \xi _{k}}\sigma _{1\lambda }(\xi
)\right\Vert _{B(E)} & \leq & \displaystyle C\left\Vert (A+(\lambda +L(\xi
))^{-1}\right\Vert _{B(E)}\left\Vert A[A+(\lambda +L(\xi ))]^{-1}\right\Vert
_{B(E)} \\ 
\vspace{-3mm} &  &  \\ 
& \leq & \displaystyle C\left[ 1+\sum\limits_{k=1}^{N}|\xi _{k}|^{2l}\right]
(1+|\lambda |+|L(\xi )|)^{-1} \\ 
\vspace{-1mm} &  &  \\ 
& \leq & \displaystyle C\left[ 1+\sum\limits_{k=1}^{N}|\xi _{k}|^{2l}\right]
\left( 1+|\lambda |+\sum\limits_{k=1}^{N}|\xi _{k}|^{2l}\right) ^{-1}~\leq
~C.%
\end{array}%
\end{equation*}%
Now let us consider the case $|\beta |=1$ and $k=1.$ Similarly, we have 
\begin{equation*}
\begin{array}{lll}
\displaystyle\left\Vert |\xi |\frac{\partial }{\partial \xi _{k}}\sigma
_{1\lambda }(\xi )\right\Vert _{B(E)} & \!\leq \! & 
\begin{array}{l}
\displaystyle C|\xi |\sum\limits_{|\alpha |\leq 2l}|\xi _{1}|^{\alpha
_{1}}\cdots |\xi _{k}|^{\alpha _{k}-1}\cdots |\xi _{N}|^{\alpha _{N}} \\ 
\times \left\Vert A[A+(\lambda +L(\xi ))]^{-2}\right\Vert _{B(E)}%
\end{array}
\\ 
\vspace{-3mm} &  &  \\ 
& \!\leq \! & \displaystyle C\sum\limits_{j=1}^{N}|\xi
_{j}|\sum\limits_{|\alpha |\leq 2l}|\xi _{1}|^{\alpha _{1}}\cdots |\xi
_{k}|^{\alpha _{k}-1}\cdots |\xi _{N}|^{\alpha _{N}}\left\Vert A[A+(\lambda
+L(\xi ))]^{-2}\right\Vert _{B(E)} \\ 
\vspace{-3mm} &  &  \\ 
& \!\leq \! & \displaystyle C\sum\limits_{|\alpha |\leq
2l}\sum\limits_{j=1}^{N}|\xi _{j}||\xi _{1}|^{\alpha _{1}}\cdots |\xi
_{k}|^{\alpha _{k}-1}\cdots |\xi _{N}|^{\alpha _{N}}\left\Vert [A+(\lambda
+L(\xi ))]^{-1}\right\Vert _{B(E)} \\ 
\vspace{-3mm} &  &  \\ 
& \!\leq \! & \displaystyle C\left[ 1+\sum\limits_{k=1}^{N}|\xi _{k}|^{2l}%
\right] \left( 1+|\lambda |+\sum\limits_{k=1}^{N}|\xi _{k}|^{2l}\right)
^{-1}~\leq ~C.%
\end{array}%
\end{equation*}%
Hence, 
\begin{equation*}
\left\Vert |\xi |\frac{\partial }{\partial \xi _{k}}\sigma _{1\lambda }(\xi
)\right\Vert _{B(E)}~\leq ~C\left[ 1+\sum\limits_{k=1}^{N}|\xi _{k}|^{2l}%
\right] \left( 1+|\lambda |+\sum\limits_{k=1}^{N}|\xi _{k}|^{2l}\right)
^{-1}~\leq ~C.
\end{equation*}%
Other cases can be proved analogously. Therefore, we obtain 
\begin{equation*}
\left\Vert |\xi |^{|\beta |}D^{\beta }\sigma _{1\lambda }(\xi )\right\Vert
_{L_{\infty }(B(E)))}~\leq ~C\eqno(30)
\end{equation*}%
and 
\begin{equation*}
\left\Vert |\xi |^{|\beta |}D^{\beta }\sigma _{2\lambda }(\xi )\right\Vert
_{L_{\infty }(B(E)))}~\leq ~K.\eqno(31)
\end{equation*}%
for each multi--index $\beta $ with $|\beta |\leq \left\lceil \frac{N}{\eta }%
\right\rceil +1\leq N+1.$ Taking into account assumptions of the theorem and
using (30), (31) and Corollary 2.13, we find $\sigma _{j\lambda }(\xi )$ are
Fourier multipliers from $B_{q_{1},r}^{s}(R^{N};E)$ to $%
B_{q_{2},r}^{s}(R^{N};E).$ Hence, for all $f\in B_{q_{1},r}^{s}(R^{N};E)$
there is a unique solution of the equation (20) in the form $%
u(x)=F^{-1}[A+(\lambda +L(\xi ))]^{-1}\hat{f}$ and the estimate (21) holds.

\hbox{\vrule height7pt width5pt}

Let $B$ denote the space $%
B(B_{q_{1},r}^{s}(R^{N};E),B_{q_{2},r}^{s}(R^{N};E))$ and $Q$ be an operator
generated by the problem (21), i.e, 
\begin{equation*}
D(Q)~=~B_{q_{2},r}^{2l,s}(R^{N};E(A),E),~~~Qu~=~\sum\limits_{|\alpha |\leq
2l}a_{\alpha }D^{\alpha }u+Au.\eqno(32)
\end{equation*}

\vspace{3mm}

\textbf{Result 4.2.} Assume all conditions of the Theorem 4.1 hold. Then for
all $\lambda \in S(\varphi )$ the resolvent of operator $Q$ exist and the
following estimate holds 
\begin{equation*}
\sum\limits_{|\alpha |\leq 2l}|\lambda |^{1-\frac{|\alpha |}{2l}}\left\Vert
D^{\alpha }(Q+\lambda )^{-1}\right\Vert _{B}+\left\Vert A(Q+\lambda
)^{-1}\right\Vert _{B}\leq C.
\end{equation*}%
\hfill

\vspace{3mm}

\textbf{Remark 4.3. }The Result 4.2 particularly, implies that the operator $%
Q+a $, $a>0$ is positive ($B_{q_{2},r}^{2l,s}\left( R;E\left( A\right)
,E\right) \rightarrow B_{q_{1},r}^{s}(R;E)$)$.$ I.e. if $A$ is strongly
positive for $\varphi \in \left( \frac{\pi }{2},\pi \right) $ then ( see
e.g. $\left[ \text{3, Theorem 8.8}\right] $ $)$ the operator $Q+a$ is a
generator of analytic semigroup$.$

\section*{5. Convolution operator equations}

In this section we shall investigate an ordinary convolution operator
equation 
\begin{equation*}
\ \left( L+\lambda \right) u=\sum\limits_{k=0}^{l}a_{k}\ast \frac{d^{k}u}{%
dx^{k}}+A_{\lambda }\ast u=f(t)\eqno(33)
\end{equation*}%
in $B_{q_{1},r}^{s}(R;E),$ where $A_{\lambda }=A_{\lambda }(x)=A(x)+\lambda $
is a possible unbounded operator in $E$ and $a_{k}=a_{k}(x)$ are complex
valued functions.

\textbf{Condition 5.1. }Suppose 
\begin{equation*}
a_{k}\in L_{1}(R)\text{, }L\left( \xi \right) =\sum\limits_{k=0}^{l}\hat{a}%
_{k}(\xi )\left( i\xi \right) ^{k}\in S\left( \varphi _{1}\right) \text{, }%
\varphi _{1}+\varphi <\pi .
\end{equation*}%
and there is a positive constant $C$ so that \ 
\begin{equation*}
\left\vert L\left( \xi \right) \right\vert >C\left\vert \xi \right\vert
^{l}\sum\limits_{k=0}^{l}\left\vert \hat{a}_{k}\right\vert .
\end{equation*}%
\textbf{\ }

\vspace{3mm}

\textbf{Lemma 5.2. }Let the Condition 5.1 be satisfied and $A(\xi )$ be a
uniformly $\varphi $-positive ($\varphi \in \lbrack 0,\pi )$) operator in a
Banach space $E,$ $\lambda \in S\left( \varphi \right) $. Then, operator
functions 
\begin{eqnarray*}
\sigma _{0}\left( \xi ,\lambda \right) &=&\lambda \left[ A\left( \xi \right)
+(\lambda +L\left( \xi \right) \right] ^{-1}\text{, } \\
\sigma _{1}\left( \xi ,\lambda \right) &=&A\left( \xi \right) \left[ A\left(
\xi \right) +(\lambda +L\left( \xi \right) \right] ^{-1}
\end{eqnarray*}

\textbf{\ } 
\begin{equation*}
\sigma _{2}\left( \xi ,\lambda \right) =\sum\limits_{k=0}^{l}\left\vert
\lambda \right\vert ^{1-\frac{k}{l}}\hat{a}_{k}(\xi )\left( i\xi \right) ^{k}%
\left[ A\left( \xi \right) +(\lambda +L\left( \xi \right) )\right] ^{-1}
\end{equation*}%
are uniformly bounded. \vspace{1mm} 

\textbf{Proof. }Let us note that for the sake of simplicity we shall not
change constants in every step.\textbf{\ }By using the resolvent properties
of positive operators we obtain 
\begin{equation*}
\left\Vert \sigma _{0}\left( \xi ,\lambda \right) \right\Vert _{B\left(
E\right) }\leq M\left\vert \lambda \right\vert (1+\left\vert \lambda
+L\left( \xi \right) \right\vert )^{-1}\leq M,
\end{equation*}%
\textbf{\ } 
\begin{eqnarray*}
\left\Vert \sigma _{1}\left( \xi ,\lambda \right) \right\Vert _{B\left(
E\right) } &=&\left\Vert A\left( \xi \right) \left[ A\left( \xi \right)
+(\lambda +L\left( \xi \right) )\right] ^{-1}\right\Vert _{B(E)} \\
&=&\left\Vert I-\left( \lambda +L\left( \xi \right) \right) \left[ A\left(
\xi \right) +(\lambda +L\left( \xi \right) )\right] ^{-1}\right\Vert _{B(E)}
\\
&\leq &1+\left\vert \lambda +L\left( \xi \right) \right\vert \left\Vert 
\left[ A\left( \xi \right) +(\lambda +L\left( \xi \right) )\right]
^{-1}\right\Vert _{B(E)} \\
&\leq &1+M\left\vert \lambda +L\left( \xi \right) \right\vert (1+\left\vert
\lambda +L\left( \xi \right) \right\vert )^{-1}\leq 1+M.
\end{eqnarray*}

Next, let us consider $\sigma _{2}.$ It is clear to see that%
\begin{equation*}
\begin{array}{lll}
\left\Vert \sigma _{2}\left( \xi \right) \right\Vert _{B\left( E\right) } & =
& \displaystyle\left\Vert \sum\limits_{k=0}^{l}\left\vert \lambda
\right\vert ^{1-\frac{k}{l}}\hat{a}_{k}(\xi )\left( i\xi \right) ^{k}\left[
A\left( \xi \right) +(\lambda +L\left( \xi \right) )\right] ^{-1}\right\Vert
_{B(E)} \\ 
\vspace{-3mm} &  &  \\ 
& \leq  & \displaystyle\leq C\sum\limits_{k=0}^{l}\left\vert \hat{a}_{k}(\xi
)\right\vert \left\vert \lambda \right\vert \left[ \left\vert \xi
\right\vert \left\vert \lambda \right\vert ^{-\frac{1}{l}}\right]
^{k}\left\Vert \left[ A\left( \xi \right) +(\lambda +L\left( \xi \right) )%
\right] ^{-1}\right\Vert _{B\left( E\right) }.%
\end{array}%
\end{equation*}

Therefore, $\sigma _{2}\left( \xi ,\lambda \right) $ is bounded if 
\begin{equation*}
\left\Vert I\right\Vert _{B\left( E\right) }=\sum\limits_{k=0}^{l}\left\vert 
\hat{a}_{k}(\xi )\right\vert \left\vert \lambda \right\vert \left\vert \xi
\right\vert ^{k}\left\vert \lambda \right\vert ^{-\frac{k}{l}}\left\Vert %
\left[ A\left( \xi \right) +(\lambda +L\left( \xi \right) )\right]
^{-1}\right\Vert _{B\left( E\right) }\leq C.
\end{equation*}%
Since $A$ is a uniformly $\varphi $-positive and $L(\xi )\in S(\varphi _{1})$
for all $\xi \in R$ then 
\begin{equation*}
\left\Vert I\right\Vert _{B\left( E\right) }\leq
C\sum\limits_{k=0}^{l}\left\vert \hat{a}_{k}(\xi )\right\vert \left\vert
\lambda \right\vert \left[ 1+\left\vert \xi \right\vert ^{l}\left\vert
\lambda \right\vert ^{-1}\right] \left[ 1+\left\vert \lambda +L\left( \xi
\right) \right\vert \right] ^{-1}.
\end{equation*}

Taking into account $\hat{a}_{k}(\xi )\left( i\xi \right) ^{k}\in S(\varphi
_{2}),$ $\varphi _{2}<\frac{\pi }{l},$ $\hat{a}_{l}(\xi )\neq 0$ and by
using $\left[ 7,\text{ Lemma 2.3}\right] ,$ Condition 5.1 and the fact $\hat{%
a}_{k}(\xi )\in L_{\infty }(R)$ (Housdorff-Youngs inequality) we get 
\begin{equation*}
\begin{array}{lll}
\left\Vert \sigma _{2}\left( \xi \right) \right\Vert _{B\left( E\right) } & 
\leq  & \displaystyle C\left\Vert I\right\Vert _{B\left( E\right) }\leq
C\sum\limits_{k=0}^{l}\left\vert \hat{a}_{k}(\xi )\right\vert \left[
\left\vert \lambda \right\vert +\left\vert \xi \right\vert ^{l}\right] \left[
1+\left\vert \lambda \right\vert +\left\vert L\left( \xi \right) \right\vert %
\right] ^{-1} \\ 
\vspace{-3mm} &  &  \\ 
& \leq  & \displaystyle C\sum\limits_{k=0}^{l}\left\vert \hat{a}_{k}(\xi
)\right\vert \left[ \left\vert \lambda \right\vert +\left\vert \xi
\right\vert ^{l}\right] \left[ 1+\left\vert \lambda \right\vert
+\sum\limits_{k=0}^{l}\left\vert \hat{a}_{k}(\xi )\right\vert \left\vert \xi
\right\vert ^{k}\right] ^{-1}\leq C.%
\end{array}%
\eqno(34)
\end{equation*}

\vspace{3mm}

\textbf{Proposition 5.3. }Let the Condition 5.1 be satisfied and $A(\xi )$
be a uniformly $\varphi $-positive ($\varphi \in \lbrack 0,\pi )$) operator
in a Banach space $E$ and 
\begin{equation*}
a_{k}\in C^{\left( m\right) }\left( R\right) ,\text{ }k=0,1,...,l\text{, }%
m=1,2,\text{ }A\left( \xi \right) A^{-1}\left( \xi _{0}\right) \in C^{\left(
m\right) }\left( R;B\left( E\right) \right) ,\text{ }\xi _{0}\in R.
\end{equation*}%
Suppose there are positive constants $C_{i}$, $i=1,...,4$ so that 
\begin{equation*}
\left\Vert A^{\left( m\right) }\left( \xi \right) A^{-1}\left( \xi \right)
\right\Vert _{B(E)}\leq C_{1},\text{ }\left\Vert \xi ^{m}A^{\left( m\right)
}\left( \xi \right) A^{-1}\left( \xi \right) \right\Vert _{B(E)}\leq C_{2}%
\eqno(35)
\end{equation*}%
\begin{equation*}
\text{ }\left\vert \xi ^{m}\hat{a}_{k}(\xi )\right\vert \leq M\text{, }%
\left\vert \frac{d^{m}}{d\xi ^{m}}\hat{a}_{k}(\xi )\right\vert \leq
C_{3},\left\vert \xi ^{m}\frac{d^{m}}{d\xi ^{m}}\hat{a}_{k}(\xi )\right\vert
\leq C_{4},\eqno(36)\text{ }
\end{equation*}%
Then, operator functions $\frac{d^{m}}{d\xi ^{m}}\sigma _{i}\left( \xi
\right) $, $i=0,1,2$ are uniformly bounded for $\lambda \in S_{\varphi }$
with $0<\lambda _{0}\leq \left\vert \lambda \right\vert .$

\vspace{1mm} \textbf{Proof. }Let us first prove for $\frac{d}{d\xi }\sigma
_{1}\left( \xi \right) $. Really, 
\begin{equation*}
\left\Vert \frac{d}{d\xi }\sigma _{1}\left( \xi \right) \right\Vert
_{B(E)}\leq \left\Vert I_{1}\right\Vert _{B(E)}+\left\Vert I_{2}\right\Vert
_{B(E)}+\left\Vert I_{3}\right\Vert _{B(E)}
\end{equation*}%
where 
\begin{eqnarray*}
I_{1} &=&A^{\prime }\left( \xi \right) \left[ A\left( \xi \right) +\lambda
+L(\xi )\right] ^{-1},\text{ } \\
I_{2} &=&A\left( \xi \right) A^{\prime }\left( \xi \right) \left[ A\left(
\xi \right) +\lambda +L(\xi )\right] ^{-2}
\end{eqnarray*}%
and 
\begin{equation*}
I_{3}=A\left( \xi \right) L^{\prime }\left( \xi \right) \left[ A\left( \xi
\right) +\lambda +L(\xi )\right] ^{-2}.
\end{equation*}%
By using\ $\left( 35\right) $ we get 
\begin{eqnarray*}
\left\Vert I_{1}\right\Vert _{B(E)} &=&\left\Vert A^{\prime }\left( \xi
\right) A^{-1}\left( \xi \right) A\left( \xi \right) \left[ A\left( \xi
\right) +\lambda +L(\xi )\right] ^{-1}\right\Vert _{B(E)} \\
&\leq &\left\Vert A^{\prime }\left( \xi \right) A^{-1}\left( \xi \right)
\right\Vert _{B(E)}\left\Vert A\left( \xi \right) \left[ A\left( \xi \right)
+\lambda +L(\xi )\right] ^{-1}\right\Vert _{B(E)} \\
&\leq &C.
\end{eqnarray*}%
Taking into account the fact $A$ is closed and linear operator and by using $%
\left( 35\right) $ we obtain 
\begin{equation*}
\left\Vert I_{2}\right\Vert _{B(E)}\leq \left\Vert A^{\prime }\left( \xi
\right) \left[ A\left( \xi \right) +\lambda +L(\xi )\right] ^{-1}\right\Vert
_{B(E)}
\end{equation*}%
\begin{equation*}
\cdot \left\Vert A\left( \xi \right) \left[ A\left( \xi \right) +\lambda
+L(\xi )\right] ^{-1}\right\Vert _{B(E)}\leq C.
\end{equation*}%
Since, $A(\xi )$ is a uniformly $\varphi $-positive, $L\left( \xi \right)
\in S\left( \varphi _{1}\right) $ for $\varphi _{1}+\varphi <\pi ,$ we get 
\begin{eqnarray*}
\left\Vert I_{3}\right\Vert _{B(E)} &\leq &\left\vert L^{\prime }\left( \xi
\right) \right\vert \left\Vert \left[ A\left( \xi \right) +\lambda +L(\xi )%
\right] ^{-1}\right\Vert _{B(E)} \\
\cdot \left\Vert A\left( \xi \right) \left[ A\left( \xi \right) +\lambda
+L(\xi )\right] ^{-1}\right\Vert _{B(E)} &\leq &C\left\vert L^{\prime
}\left( \xi \right) \right\vert \left[ 1+\left\vert \lambda +L(\xi
)\right\vert \right] ^{-1}.
\end{eqnarray*}%
Thus, by using $\left[ 7,\text{ Lemma 2.3}\right] $ we have%
\begin{equation*}
\left\Vert I_{3}\right\Vert _{B(E)}\leq C\left\vert L^{\prime }\left( \xi
\right) \right\vert \left[ 1+\left\vert \lambda \right\vert
+\sum\limits_{k=0}^{l}\left\vert \hat{a}_{k}(\xi )\right\vert \left\vert \xi
\right\vert ^{k}\right] ^{-1}.
\end{equation*}%
It is clear to see that 
\begin{equation*}
\left\vert L^{\prime }\left( \xi \right) \right\vert \leq
\sum\limits_{k=0}^{l}\left\vert \frac{d}{d\xi }\hat{a}_{k}(\xi )\right\vert
\left\vert \xi \right\vert ^{k}+l\sum\limits_{k=1}^{l}\left\vert \hat{a}%
_{k}(\xi )\right\vert \left\vert \xi \right\vert ^{k-1}.\eqno(37)
\end{equation*}

By using $\left( 36\right) $, $\left( 37\right) $ we obtain 
\begin{equation*}
\sum\limits_{k=0}^{l}\left\vert \frac{d}{d\xi }\hat{a}_{k}(\xi )\right\vert
\left\vert \xi \right\vert ^{k}\leq C\left[ 1+\left\vert \lambda \right\vert
+\sum\limits_{k=0}^{l}\left\vert \hat{a}_{k}(\xi )\right\vert \left\vert \xi
\right\vert ^{k}\right]
\end{equation*}%
and 
\begin{equation*}
\sum\limits_{k=1}^{l}\left\vert \hat{a}_{k}(\xi )\right\vert \left\vert \xi
\right\vert ^{k-1}\leq C\left[ 1+\left\vert \lambda \right\vert
+\sum\limits_{k=0}^{l}\left\vert \hat{a}_{k}(\xi )\right\vert \left\vert \xi
\right\vert ^{k}\right]
\end{equation*}

that implies 
\begin{equation*}
\left\Vert I_{3}\right\Vert _{B(E)}\leq C\left\vert L^{\prime }\left( \xi
\right) \right\vert \left[ 1+\left\vert \lambda \right\vert
+\sum\limits_{k=0}^{l}\left\vert \hat{a}_{k}(\xi )\right\vert \left\vert \xi
\right\vert ^{k}\right] ^{-1}\leq C.\eqno(38)
\end{equation*}%
Next we shall prove uniformly boundedness of $\frac{d}{d\xi }\sigma
_{2}\left( \xi ,\lambda \right) $. Similarly, 
\begin{equation*}
\left\Vert \frac{d}{d\xi }\sigma _{2}\left( \xi \right) \right\Vert
_{B(E)}\leq \left\Vert J_{1}\right\Vert _{B(E)}+\left\Vert J_{2}\right\Vert
_{B(E)}+\left\Vert J_{3}\right\Vert _{B(E)}+\left\Vert J_{4}\right\Vert
_{B(E)},
\end{equation*}%
where 
\begin{eqnarray*}
J_{1} &=&\sum\limits_{k=0}^{l}\left\vert \lambda \right\vert ^{1-\frac{k}{l}}%
\frac{d}{d\xi }\hat{a}_{k}(\xi )\left( i\xi \right) ^{k}\left[ A\left( \xi
\right) +L(\xi )\right] ^{-1},\text{ } \\
J_{2} &=&\sum\limits_{k=0}^{l}\left\vert \lambda \right\vert ^{1-\frac{k}{l}}%
\hat{a}_{k}(\xi )ik\left( i\xi \right) ^{k-1}\left[ A\left( \xi \right)
+L(\xi )\right] ^{-1}, \\
J_{3} &=&\sum\limits_{k=0}^{l}\left\vert \lambda \right\vert ^{1-\frac{k}{l}}%
\hat{a}_{k}(\xi )\left( i\xi \right) ^{k}L^{\prime }\left( \xi \right) \left[
A\left( \xi \right) +L(\xi )\right] ^{-2}
\end{eqnarray*}%
and 
\begin{equation*}
J_{4}=\sum\limits_{k=0}^{l}\left\vert \lambda \right\vert ^{1-\frac{k}{l}}%
\hat{a}_{k}(\xi )\left( i\xi \right) ^{k}A^{\prime }\left( \xi \right) \left[
A\left( \xi \right) +L(\xi )\right] ^{-2}.
\end{equation*}

Let us first show $J_{1}$ is uniformly bounded. Since, 
\begin{equation*}
\left\Vert J_{1}\right\Vert _{B(E)}\leq \sum\limits_{k=0}^{l}\left\vert 
\frac{d}{d\xi }\hat{a}_{k}(\xi )\right\vert \left\Vert \left\vert \lambda
\right\vert ^{1-\frac{k}{l}}\left( i\xi \right) ^{k}\left[ A\left( \xi
\right) +(\lambda +L(\xi ))\right] ^{-1}\right\Vert _{B(E)}
\end{equation*}%
by virtue of $\left( 34\right) $ and $\left( 37\right) $ we obtain $%
\left\Vert J_{1}\right\Vert _{B(E)}\leq C.$ Then, with the help of $\left(
34\right) $, $\left( 37\right) ,$ Condition 5.1 and the fact $\hat{a}%
_{k}(\xi )\in L_{\infty }(R)$ we get 
\begin{eqnarray*}
\left\Vert J_{2}\right\Vert _{B(E)} &\leq &\sum\limits_{k=1}^{l}\left\vert 
\hat{a}_{k}(\xi )\right\vert \left\Vert \left\vert \lambda \right\vert ^{1-%
\frac{k}{l}}\left( i\xi \right) ^{k-1}\left[ A\left( \xi \right) +(\lambda
+L(\xi ))\right] ^{-1}\right\Vert _{B(E)} \\
&\leq &C\sum\limits_{k=1}^{l}\left\vert \hat{a}_{k}(\xi )\right\vert
\left\vert \lambda \right\vert ^{-\frac{1}{l}}\left\vert \lambda \right\vert
\left( \left\vert \xi \right\vert \left\vert \lambda \right\vert ^{\frac{-1}{%
l}}\right) _{B(E)}^{k}\left( 1+\left\vert \lambda \right\vert
+\sum\limits_{k=0}^{l}\left\vert \hat{a}_{k}(\xi )\left( i\xi \right)
^{k}\right\vert \right) ^{-1}
\end{eqnarray*}%
\begin{eqnarray*}
&\leq &C\left\vert \lambda _{0}\right\vert ^{-\frac{1}{l}}\sum%
\limits_{k=1}^{l}\left\vert \hat{a}_{k}(\xi )\right\vert \left\vert \lambda
\right\vert \left( 1+\left\vert \lambda \right\vert ^{-1}\left\vert \xi
\right\vert ^{l}\right) \left( 1+\left\vert \lambda \right\vert +\left\vert
\xi \right\vert ^{l}\sum\limits_{k=0}^{l}\left\vert \hat{a}_{k}(\xi
)\right\vert \right) ^{-1} \\
&\leq &C_{0}\sum\limits_{k=1}^{l}\left\vert \hat{a}_{k}(\xi )\right\vert
\left( \left\vert \lambda \right\vert +\left\vert \xi \right\vert
^{l}\right) \left( 1+\left\vert \lambda \right\vert +\left\vert \xi
\right\vert ^{l}\sum\limits_{k=0}^{l}\left\vert \hat{a}_{k}(\xi )\right\vert
\right) ^{-1}\leq C_{0}.
\end{eqnarray*}%
Next, by means of $\left( 34\right) ,\left( 35\right) ,\left( 38\right) $
and the facts $\hat{a}_{k}(\xi )\in L_{\infty }(R),$ $A(\xi )$ is a
uniformly $\varphi $-positive, $L\left( \xi \right) \in S\left( \varphi
_{1}\right) $ for $\varphi _{1}+\varphi <\pi ,$ we obtain%
\begin{equation*}
\begin{array}{lll}
\left\Vert J_{3}\right\Vert _{B(E)} & \leq  & \displaystyle C\left\vert
L^{\prime }\left( \xi \right) \right\vert \left\Vert \left[ A\left( \xi
\right) +(\lambda +L(\xi ))\right] ^{-1}\right\Vert _{B(E)} \\ 
\vspace{-3mm} & \cdot  & \sum\limits_{k=0}^{l}\left\vert \hat{a}_{k}(\xi
)\right\vert \left\vert \lambda \right\vert \left\vert \xi \right\vert
^{k}\left\vert \lambda \right\vert ^{-\frac{k}{l}}\left\Vert \left[ A\left(
\xi \right) +(\lambda +L(\xi ))\right] ^{-1}\right\Vert _{B(E)} \\ 
&  &  \\ 
& \leq  & C\sum\limits_{k=0}^{l}\left\vert \hat{a}_{k}(\xi )\right\vert
\left\vert \xi \right\vert ^{k}\left[ 1+\left\vert \lambda \right\vert
+\sum\limits_{k=0}^{l}\left\vert \hat{a}_{k}(\xi )\right\vert \left\vert \xi
\right\vert ^{k}\right] ^{-1} \\ 
&  &  \\ 
& \leq  & C.%
\end{array}%
\eqno(39)
\end{equation*}%
Finally, by virtue of $\left( 34\right) $ and $\left( 35\right) $ we obtain 
\begin{eqnarray*}
\left\Vert J_{4}\right\Vert _{B(E)} &\leq &C\left\Vert \frac{d}{d\xi }%
A\left( \xi \right) A^{-1}\left( \xi \right) A\left( \xi \right) \left[
A\left( \xi \right) +(\lambda +L(\xi ))\right] ^{-1}\right\Vert _{B(E)} \\
&&\cdot \sum\limits_{k=0}^{l}\left\vert \hat{a}_{k}(\xi )\right\vert
\left\vert \lambda \right\vert \left\vert \xi \right\vert ^{k}\left\vert
\lambda \right\vert ^{\frac{-k}{l}}\left\Vert \left[ A\left( \xi \right)
+(\lambda +L(\xi ))\right] ^{-1}\right\Vert _{B(E)} \\
&\leq &C\left\Vert A^{\prime }\left( \xi \right) A^{-1}\left( \xi \right)
\right\Vert _{B(E)}\left\Vert A\left( \xi \right) \left[ A\left( \xi \right)
+(\lambda +L(\xi ))\right] ^{-1}\right\Vert _{B(E)} \\
&\leq &C.
\end{eqnarray*}%
Hence, operator functions $\frac{d}{d\xi }\sigma _{i}\left( \xi ,\lambda
\right) $ are uniformly bounded\textbf{. }Hence, operator functions $\frac{d%
}{d\xi }\sigma _{i}\left( \xi \right) $ are uniformly bounded\textbf{. }%
Using the same techniques one can easily establish boundedness of $\frac{%
d^{2}}{d\xi ^{2}}\sigma _{i}\left( \xi \right) .$

\vspace{3mm}

\textbf{Proposition 5.4. }Let all conditions of Proposition 5.3 are satisfied%
$.$

Then the following estimates hold 
\begin{equation*}
\left\Vert \left\vert \xi \right\vert ^{m}\frac{d^{m}}{d\xi ^{m}}\sigma
_{i}\left( \xi ,\lambda \right) \right\Vert _{L_{\infty }(B(E)))}\leq A_{i}%
\text{, }m,i=0,1,2.
\end{equation*}

\vspace{1mm}

\textbf{Proof. }As a matter of fact, it is enough to prove 
\begin{equation*}
\left\vert \xi \right\vert \left\Vert I_{i}\right\Vert _{B(E)}\leq C_{i}%
\text{ and }\left\vert \xi \right\vert \left\Vert J_{j}\right\Vert
_{B(E)}\leq D_{j}
\end{equation*}%
for some constant $C_{i}$ and $D_{j}$, $i=1,$ $2,$ $3$, $j=1,$ $2,$ $3,$ $4.$
It is easy to see from the proof of Proposition 5.3 
\begin{eqnarray*}
\left\vert \xi \right\vert \left\Vert I_{1}\right\Vert _{B(E)} &\leq
&C_{1}\left\Vert \xi A^{\prime }\left( \xi \right) A^{-1}\left( \xi \right)
\right\Vert _{B(E)} \\
&&\cdot \left\Vert A\left( \xi \right) \left[ A\left( \xi \right) +(\lambda
+L(\xi ))\right] ^{-1}\right\Vert _{B(E)}^{2} \\
\left\vert \xi \right\vert \left\Vert I_{2}\right\Vert _{B(E)} &\leq
&C_{2}\left\Vert \xi A^{\prime }\left( \xi \right) A^{-1}\left( \xi \right)
\right\Vert _{B(E)} \\
&&\cdot \left\Vert A\left( \xi \right) \left[ A\left( \xi \right) +(\lambda
+L(\xi ))\right] ^{-1}\right\Vert _{B(E)}
\end{eqnarray*}%
\begin{equation*}
\left\vert \xi \right\vert \left\Vert I_{3}\right\Vert _{B(E)}\leq
C_{3}\left\vert \xi \right\vert \left\vert L^{\prime }\left( \xi \right)
\right\vert \left[ 1+\left\vert \lambda \right\vert
+\sum\limits_{k=0}^{l}\left\vert \hat{a}_{k}(\xi )\right\vert \left\vert \xi
\right\vert ^{k}\right] ^{-1}.
\end{equation*}%
From resolvent properties of positive operators, it follows $\xi I_{1}$ and $%
\xi I_{2}$ are uniformly bounded. By using $\left( 37\right) $ and $\left(
38\right) $ we obtain 
\begin{equation*}
\left\vert \xi \right\vert \left\Vert I_{3}\right\Vert _{B(E)}\leq
C_{3}\sum\limits_{k=0}^{l}\left\vert \hat{a}_{k}(\xi )\right\vert \left\vert
\xi \right\vert ^{k}\left[ 1+\left\vert \lambda \right\vert
+\sum\limits_{k=0}^{l}\left\vert \hat{a}_{k}(\xi )\right\vert \left\vert \xi
\right\vert ^{k}\right] ^{-1}\leq C_{3}.
\end{equation*}%
Similarly, from the proof of Lemma 5.3 it follows 
\begin{equation*}
\left\vert \xi \right\vert \left\Vert J_{1}\right\Vert _{B(E)}\leq
\sum\limits_{k=0}^{l}\left\vert \xi \frac{d}{d\xi }\hat{a}_{k}(\xi
)\right\vert \left\Vert \left\vert \lambda \right\vert ^{1-\frac{k}{l}%
}\left( i\xi \right) ^{k}\left[ A\left( \xi \right) +(\lambda +L(\xi ))%
\right] ^{-1}\right\Vert _{B(E)},
\end{equation*}%
\begin{equation*}
\left\vert \xi \right\vert \left\Vert J_{2}\right\Vert _{B(E)}\leq
\sum\limits_{k=0}^{l}\left\vert \hat{a}_{k}(\xi )\right\vert \left\Vert
\left\vert \lambda \right\vert ^{1-\frac{k}{l}}\left( i\xi \right) ^{k}\left[
A\left( \xi \right) +(\lambda +L(\xi ))\right] ^{-1}\right\Vert _{B(E)},
\end{equation*}%
\begin{eqnarray*}
\left\vert \xi \right\vert \left\Vert J_{3}\right\Vert _{B(E)} &\leq
&C\left\vert \xi L^{\prime }\left( \xi \right) \right\vert \left\Vert \left[
A\left( \xi \right) +(\lambda +L(\xi ))\right] ^{-1}\right\Vert _{B(E)} \\
&&\cdot \sum\limits_{k=0}^{l}\left\vert \left\vert \hat{a}_{k}(\xi
)\right\vert \lambda \right\vert \left\vert \xi \right\vert ^{k}\left\vert
\lambda \right\vert ^{\frac{-k}{l}}\left\Vert \left[ A\left( \xi \right)
+(\lambda +L(\xi ))\right] ^{-1}\right\Vert _{B(E)},
\end{eqnarray*}%
and 
\begin{equation*}
\left\vert \xi \right\vert \left\Vert J_{4}\right\Vert _{B(E)}\leq
C\left\Vert \xi A^{\prime }\left( \xi \right) A^{-1}\left( \xi \right)
\right\Vert _{B(E)}\left\Vert A\left( \xi \right) \left[ A\left( \xi \right)
+(\lambda +L(\xi ))\right] ^{-1}\right\Vert _{B(E)}.
\end{equation*}

Using $\left( 34\right) $, $\left( 35\right) $, $\left( 36\right) $, $\left(
38\right) $ and the fact $\hat{a}_{k}(\xi )\in L_{\infty }(R)$ it is easy to
show that $\left\vert \xi \right\vert \left\Vert J_{1}\right\Vert _{B(E)}$, $%
\left\vert \xi \right\vert \left\Vert J_{2}\right\Vert _{B(E)},$ $\left\vert
\xi \right\vert \left\Vert J_{3}\right\Vert _{B(E)}$ and $\left\vert \xi
\right\vert \left\Vert J_{4}\right\Vert _{B(E)}$ are uniformly bounded$.$ By
virtue of the same techniques one can easily establish uniformly boundedness
of $\left\vert \xi \right\vert ^{m}\frac{d^{2}}{d\xi ^{2}}\sigma _{i}\left(
\xi \right) \ $for $m=1,2.$ \ \hbox{\vrule height7pt width5pt}

\vspace{3mm}

\textbf{Corollary 5.5. }Assume all conditions of Proposition 5.4 are
satisfied. Then, operator-functions $\sigma _{i}\left( \xi \right) $ are
Fourier multipliers from $B_{q_{1},r}^{s}(R^{n};E)$ to $%
B_{q_{2},r}^{s}(R^{n};E).$

\vspace{1mm}

\textbf{Proof. }To prove $\sigma _{i}\left( \xi \right) $ are uniformly
bounded multipliers from $B_{q_{1},r}^{s}(R^{n};E)$ to $%
B_{q_{2},r}^{s}(R^{n};E),$ we need to show $\sigma _{i}\in C^{\left(
1\right) }(R;B(E))$ and there exists a constant \ $K>0$ such that, 
\begin{equation*}
\left\Vert (1+\left\vert \xi \right\vert )^{\left\vert \beta \right\vert }%
\frac{d^{\beta }}{d\xi ^{\beta }}\sigma _{i}\left( \xi \right) \right\Vert
_{L_{\infty }(R;B(E))}\leq K
\end{equation*}%
for each multi--index $\beta $ with $|\beta |\leq \left\lceil \frac{1}{\eta }%
\right\rceil +1\leq 2.$ From the Proposition 5.2, Proposition 5.3 and
Proposition 5.4 it follows $\sigma _{i}\in C^{1}(R;B(E))$ and 
\begin{equation*}
\left\Vert \frac{d^{m}}{d\xi ^{m}}\sigma _{i}\left( \xi \right) \right\Vert
_{L_{\infty }(B(E)))}\leq A_{1},\left\Vert \left\vert \xi \right\vert ^{m}%
\frac{d^{m}}{d\xi ^{m}}\sigma _{i}\left( \xi \right) \right\Vert _{L_{\infty
}(B(E)))}\leq A_{2}
\end{equation*}%
for every $i,m=0,1,2.$ Hence, $\sigma _{i}\left( \xi \right) $ are Fourier
multipliers from $B_{q_{1},r}^{s}(R^{n};E)$ to $B_{q_{2},r}^{s}(R^{n};E).$

\vspace{3mm}

\hbox{\vrule height7pt width5pt}

\vspace{3mm}

\textbf{Theorem 5.6.}\ \textbf{\ }Let $\ f\in B_{q_{1},r}^{s}(R;E)$ and$%
\frac{1}{q_{2}}=\frac{1}{q_{1}}-\frac{1}{\eta ^{\prime }},$ $1<q_{1}<\eta
^{\prime }\leq \infty $ $.$ Then, $\left( 33\right) $ has a unique solution $%
u\in B_{q_{2},r}^{l,s}\left( R;E\left( A\right) ,E\right) $ and the
following coercive uniform estimate holds%
\begin{equation*}
\begin{array}{lll}
\left\Vert \lambda u\right\Vert _{B_{q_{2},r}^{s}\left( R;E\left( A\right)
,E\right) } & + & \displaystyle\sum\limits_{k=0}^{l}\left\vert \lambda
\right\vert ^{1-\frac{k}{l}}\left\Vert a_{k}\ast \frac{d^{k}u}{dx^{k}}%
\right\Vert _{B_{q_{2},r}^{s}(R;E)} \\ 
\vspace{-3mm} & + & \left\Vert A\ast u\right\Vert
_{B_{q_{2},r}^{s}(R;E)}\leq C\left\Vert f\right\Vert _{B_{q_{1},r}^{s}(R;E)}%
\end{array}%
\eqno(40)
\end{equation*}

for all $\lambda \in S_{\varphi }$ with $0<\lambda _{0}\leq \left\vert
\lambda \right\vert $, provided the bellow conditions satisfied:

(1) $E$ is a Banach space$;$

(2) Condition 5.1 holds and 
\begin{equation*}
\hat{a}_{k}\in C^{\left( m\right) }\left( R\right) ,\text{ }k=0,1,...,l\text{%
, }\hat{A}\left( \xi \right) \hat{A}^{-1}\left( \xi _{0}\right) \in
C^{\left( m\right) }\left( R;B\left( E\right) \right) \text{, }\xi _{0}\in R;
\end{equation*}

(3) $\hat{A}(\xi )$ is a uniformly $\varphi $-positive ($\varphi \in \lbrack
0,\pi )$) operator in $E$. Moreover, there are positive constants $C_{i}$, $%
i=1,...,4$ so that for $m=0,1,2$ 
\begin{equation*}
\text{ }\left\vert \xi ^{m}\hat{a}_{k}(\xi )\right\vert \leq M\text{, }%
\left\vert \frac{d^{m}}{d\xi ^{m}}\hat{a}_{k}(\xi )\right\vert \leq
C_{1},\left\vert \xi ^{m}\frac{d^{m}}{d\xi ^{m}}\hat{a}_{k}(\xi )\right\vert
\leq C_{2};\text{ }
\end{equation*}

\begin{equation*}
\left\Vert \hat{A}^{\left( m\right) }\left( \xi \right) \hat{A}^{-1}\left(
\xi \right) \right\Vert _{B(E)}\leq C_{3},\text{ }\left\Vert \xi ^{m}\hat{A}%
^{\left( m\right) }\left( \xi \right) \hat{A}^{-1}\left( \xi \right)
\right\Vert _{B(E)}\leq C_{4}.
\end{equation*}

\ \textbf{Proof}. Applying Fourier transform to equation $\left( 37\right) $%
\ we get 
\begin{equation*}
\left[ \hat{A}\left( \xi \right) +L\left( \xi \right) \right] \hat{u}\left(
\xi \right) =\hat{f}\left( \xi \right) .
\end{equation*}

Since\ $L(\xi )\in S(\varphi _{1})$ for all $\xi \in R$ and $\hat{A}$ is
positive$,$ the operator\ $\hat{A}\left( \xi \right) +L\left( \xi \right) $
is invertible in $E$. Thus, we obtain that the solution of equation $\left(
33\right) $ can be represented in the following form 
\begin{equation*}
u\left( x\right) =F^{-1}\left[ \hat{A}\left( \xi \right) +\lambda +L\left(
\xi \right) \right] ^{-1}\hat{f}.\eqno(41)
\end{equation*}

By using $\left( 41\right) $ we get

\begin{eqnarray*}
\left\Vert A\ast u\right\Vert _{B_{q_{2},r}^{s}(R;E)} &=&\left\Vert F^{-1} 
\left[ \sigma _{1}\left( \xi \right) \hat{f}\right] \right\Vert
_{B_{q_{2},r}^{s}(R;E)} \\
\sum\limits_{k=0}^{l}\left\vert \lambda \right\vert ^{1-\frac{k}{l}%
}\left\Vert a_{k}\ast \frac{d^{k}u}{dx^{k}}\right\Vert
_{B_{q_{2},r}^{s}(R;E)} &=&\left\Vert F^{-1}\left[ \sigma _{2}\left( \xi
\right) \hat{f}\right] \right\Vert _{B_{q_{2},r}^{s}(R;E)},
\end{eqnarray*}%
where 
\begin{eqnarray*}
\sigma _{0}\left( \xi \right) &=&\left[ \hat{A}\left( \xi \right) +\lambda
+L\left( \xi \right) \right] ^{-1},\text{ }\sigma _{1}\left( \xi \right) =%
\hat{A}\left( \xi \right) \left[ \hat{A}\left( \xi \right) +\lambda +L\left(
\xi \right) \right] ^{-1},\text{ } \\
\sigma _{2}\left( \xi \right) &=&\sum\limits_{k=0}^{l}\left\vert \lambda
\right\vert ^{1-\frac{k}{l}}\hat{a}_{k}(\xi )\left( i\xi \right) ^{k}\left[ 
\hat{A}\left( \xi \right) +\lambda +L\left( \xi \right) \right] ^{-1}.
\end{eqnarray*}

From Corollary 5.5 we know that operator-functions $\sigma _{i}\left( \xi
\right) $ are uniformly bounded multipliers from $B_{q_{1},r}^{s}(R^{n};E)$
to $B_{q_{2},r}^{s}(R^{n};E).$

Since, 
\begin{equation*}
\left\Vert A\ast u\right\Vert _{B_{q_{2},r}^{s}(R;E)}\leq C_{1}\left\Vert
f\right\Vert _{B_{q_{1},r}^{s}(R;E)},
\end{equation*}%
\begin{equation*}
\sum\limits_{k=0}^{l}\left\vert \lambda \right\vert ^{1-\frac{k}{l}%
}\left\Vert a_{k}\ast \frac{d^{k}u}{dx^{k}}\right\Vert
_{B_{q_{2},r}^{s}(R;E)}\leq C_{2}\left\Vert f\right\Vert
_{B_{q_{1},r}^{s}(R;E)}
\end{equation*}%
we obtain that there is a unique solution of the equation $\left( 33\right) $
in the form $u\left( x\right) =F^{-1}\left[ A+\lambda +L\left( \xi \right) %
\right] ^{-1}\hat{f}$ and the estimate $\left( 40\right) $\ holds for all \ $%
f\in B_{q_{1},r}^{s}(R;E).$

\hbox{\vrule height7pt width5pt}

\vspace{3mm}

Let $Q$ be an operator that generates the problem $\left( 33\right) $ i. e. 
\begin{equation*}
D\left( Q\right) =B_{q_{2},r}^{l,s}\left( R;E\left( A\right) ,E\right) ,%
\text{ }Qu=\sum\limits_{k=0}^{l}a_{k}\ast \frac{d^{k}u}{dx^{k}}+A_{\lambda
}\ast u.
\end{equation*}%
\hbox{\vrule height7pt width5pt}

\textbf{Result 5.7. }Assume all conditions of the Theorem 5.6 hold. Then for
all $\lambda \in S_{\varphi }$ with $0<\lambda _{0}\leq \left\vert \lambda
\right\vert $ the resolvent of operator $Q$ exist and the following estimate
holds 
\begin{equation*}
\sum\limits_{k=0}^{l}|\lambda |^{1-\frac{k}{l}}\left\Vert a_{k}\ast \left[ 
\frac{d^{k}}{dx^{k}}\left( Q+\lambda \right) ^{-1}\right] \right\Vert
_{B(B_{q_{1},r}^{s},\text{ }B_{q_{2},r}^{s})}+
\end{equation*}%
\begin{equation*}
\left\Vert \left\vert \lambda \right\vert \left( Q+\lambda \right)
^{-1}\right\Vert _{B(B_{q_{1},r}^{s},\text{ }B_{q_{2},r}^{s})}+\left\Vert
A\ast \left( Q+\lambda \right) ^{-1}\right\Vert _{B(B_{q_{1},r}^{s},\text{ }%
B_{q_{2},r}^{s})}\leq C.
\end{equation*}

\vspace{3mm}

\textbf{Remark 5.8. }The Result 5.7 particularly, implies that the operator $%
Q+a $, $a>0$ is positive ($B_{q_{2},r}^{s}\left( R;E\right) \rightarrow
B_{q_{1},r}^{s}(R;E)$)$.$ I.e. if $A$ is strongly positive for $\varphi \in
\left( \frac{\pi }{2},\pi \right) $ then ( see e.g. $\left[ \text{3, Theorem
8.8}\right] $ $)$ the operator $Q+a$ is a generator of analytic semigroup$.$

\section*{6. Infinite systems of quasi elliptic equations}

Consider the following infinite system 
\begin{equation*}
\sum\limits_{|\alpha |\leq 2l}a_{\alpha }D^{\alpha
}u_{m}+\sum\limits_{j=1}^{\infty }\left( d_{j}+\lambda \right)
u_{j}(x)~=~f_{m}(x),~~~x\in R^{N},~~~m=1,2,\cdots ,\infty .\eqno(41)
\end{equation*}

\vspace{1mm}

Let 
\begin{equation*}
\begin{array}{l}
D=\displaystyle\{d_m\},~~~d_m>0,~~~u=\{u_m\},~~~Du=\{d_mu_m\},~~~m=1,2,%
\cdots\infty, \\ 
\vspace{-3mm} \\ 
l_{q}(D)=\displaystyle\left\{u:~u\in
l_q,~\|u\|_{l_q(D)}=\|Du\|_{l_q}=\left(\sum\limits_{m=1}^\infty|d_mu_m|^q%
\right)^{\frac{1}{q}}<\infty \right\},~~~1<q<\infty.%
\end{array}%
\end{equation*}

\vspace{1mm}

Let $Q$ be a differential operator generating the boundary value problem
(41).

\vspace{3mm}

\textbf{Theorem 6.1.} Suppose $\frac{1}{q_{2}}=\frac{1}{q_{1}}-\frac{1}{\eta
^{\prime }},$ $1<q_{1}<\eta ^{\prime }\leq \infty $ and the following
conditions hold: 
\begin{equation*}
L(\xi )~=~\sum\limits_{|\alpha |\leq 2l}a_{\alpha }(i\xi _{1})^{\alpha
_{1}}(i\xi _{2})^{\alpha _{2}}\cdots (i\xi _{N})^{\alpha _{N}}~\in
~S(\varphi ),
\end{equation*}%
\begin{equation*}
\sum\limits_{m=1}^{\infty }d_{m}^{-1}~<~\infty ,~~~|L(\xi )|~\geq
~K\sum\limits_{k=1}^{n}|\xi _{k}|^{2l},\xi \in R^{N},~\varphi _{1}+\varphi
<\pi .
\end{equation*}

\vspace{1mm}

\noindent Then,

\vspace{2mm}

(a) For all $f(x)=\{f_{m}(x)\}_{1}^{\infty }\in
B_{q_{1},r}^{s}(R^{N};l_{q}(D))$ and $\lambda \in S(\varphi ),~\varphi \in
\lbrack 0,\pi )$ the problem (41) has a unique solution $u=\{u_{m}(x)%
\}_{1}^{\infty }$ that belongs to space $%
B_{q_{2},r}^{2l,s}(R^{N};l_{q}(D),l_{q})$ and the coercive estimate 
\begin{equation*}
\sum\limits_{|\alpha |\leq 2l}\left\Vert D^{\alpha }u\right\Vert
_{B_{q_{2},r}^{s}(R^{N};l_{q})}+\Vert Au\Vert
_{B_{q_{2},r}^{s}(R^{N};l_{q})}~\leq ~C\Vert f\Vert
_{B_{q_{1},r}^{s}(R^{N};l_{q})}\eqno(42)
\end{equation*}%
hold for the solution of the problem (41).

\vspace{2mm}

(b) There exists a resolvent $\left( Q+\lambda \right) ^{-1}$ of operator $Q$
and 
\begin{equation*}
\sum\limits_{|\alpha |\leq 2l}\left\Vert D^{\alpha }(Q+\lambda
)^{-1}\right\Vert +\left\Vert A(Q+\lambda )^{-1}\right\Vert ~\leq ~C.\eqno%
(43)
\end{equation*}

\vspace{1mm}

\textbf{Proof.} Let $E=l_{q}$ and $~A$ be an infinite matrix such that 
\begin{equation*}
A~=~[d_{m}(x)\delta _{jm}],~~~m,~j=1,2,\cdots ,\infty .
\end{equation*}%
It is clear to see that the operator $A$ is positive in $l_{q}.$ Therefore,
using Theorem 4.1 we get that the problem (41) has a unique solution $u\in
B_{q_{2},r}^{2l,s}(G;l_{q}(D),l_{q})$ and estimates (42) and (43) holds for
all $f\in B_{q_{1},r}^{s}(R^{N};l_{q})$.

\vspace{3mm}

\textbf{Remark 6.2.} There are lots of positive operators in concrete Banach
spaces. Therefore, putting concrete Banach spaces instead of $E$ and
concrete positive differential, pseudo differential operators, or finite,
infinite matrices, etc. instead of operator $A$ in (20), we can obtain the
maximal regularity of different classes of BVPs for partial differential
equations or system of equations by Theorem 4.1.


\begin{thebibliography}{99}
\bibitem{A} \textbf{H. Amann,} \textit{Linear and quasilinear parabolic
problems, Vol. I, Abstarct Linear Theory}, \textsl{Birkh\"{a}user}, Boston,
1995.

\bibitem{B} \textbf{H. Amann,} Compact embedding of vector--valued Sobolev
and Besov spaces, \textsl{Glasnik Mathematicki, \textbf{35(55) }}(2000),
161--177.

\bibitem{C} \textbf{H. Amann,} Operator--valued Fourier multipliers,
vector--valued Besov spaces, and applications, \textsl{Math. Nachr., \textbf{%
186}}(1997), 5--56.

\bibitem{D} \textbf{R.P. Agarwal, M. Bohner and V.B. Shakhmurov,} Maximal
regular boundary value problems in Banach--valued weighted spaces, \textsl{%
Boundary value problems, \textbf{1}}(2005), 9--42.

\bibitem{E} \textbf{A. Pelczynski and M. Wojciechowski,} Molecular
decompositions and embedding theorems for vector--valued Sobolev spaces with
gradient norm, \textsl{Studia Math., \textbf{107}}(1993), 61--100.

\bibitem{F} \textbf{J.P. Aubin,} Abstract boundary--value operators and
their adjoint, \textsl{Rend. Sem. Math. Univ. Padova, \textbf{43}}(1970),
1--33.

\bibitem{G} \textbf{C. Dore and S. Yakubov,} Semigroup estimates and non
coercive boundary value problems, \textsl{Semigroup Form, \textbf{60}}%
(2000), 93--121.

\bibitem{H} \textbf{R. Denk, M. Hieber and J. Pr\"{u}ss,} $R-$boundedness,
Fourier multipliers and problems of elliptic and parabolic type, \textsl{%
Mem. Amer. Math. Soc., \textbf{166}}(2003), No. 788.

\bibitem{I} \textbf{A. Favini,} Su un problema ai limiti per certa equazini
astratte del secondo ordine, \textsl{Rend. Sem. Mat. Univ. Padova, \textbf{53%
}}(1975), 211--230.

\bibitem{J} \textbf{M. Girardi and L. Weis,} Operator--valued multiplier
theorems on Besov spaces, \textsl{Math. Nachr., \textbf{251}}(2003), 34--51.

\bibitem{K} \textbf{J. Peetre,} Sur la transformation de Fourier des
fonctions a valeurs vectorielles, \textsl{Rend. Sem. Math. Univ. Padova, 
\textbf{42}}(1969),

\bibitem{L} \textbf{J. B\"{o}rgh and J. L\"{o}fstr\"{o}m,} \textit{%
Interpolation spaces: An introduction}, \textsl{Springer--Verlag}, Berlin
1976, Grundlehren der Mathematischen Wissenschaften, No. 223.

\bibitem{M} \textbf{S.G. Krein,} \textit{Linear differential equations in
Banach space}, Providence, 1982. \vspace{1mm}

\bibitem{N} \textbf{V.B. Shakhmurov,} Embedding theorems and\ maximal
regular differential operator equations in Banach-valued function spaces, 
\textsl{Inequalities and Applications., \textbf{292}}(2005), v. 2, no. 4.

\bibitem{O} \textbf{V.B. Shakhmurov} Embedding and maximal regular
differential operators in Sobolev-Lions spaces, \textsl{Acta Mathematica
Sinica., 22}(2006), 1493-1508

\bibitem{P} \textbf{T.R. McConnell,} On Fourier multiplier transformations
of Banach--valued functions, \textsl{Trans. Amer. Mat. Soc., \textbf{285}}%
(1984), 739--757.

\bibitem{Q} \textbf{H. Triebel,} \textit{Interpolation theory. Function
spaces. Differential operators}, \textsl{North--Holland}, Amsterdam, 1978.

\bibitem{R} \textbf{L. Weis,} Operator--valued Fourier multiplier theorems
and maximal $L_{p}$--regularity, \textsl{Math. Ann., }\textbf{319}(2001),
735--758.
\end{thebibliography}
\end{document}